\documentclass[reqno,date]{article}

\usepackage{amsfonts}
\usepackage{amssymb}
\usepackage{times}
\usepackage{mathptmx}
\usepackage{graphicx}
\usepackage[hmargin=1.5in, vmargin=1.5in]{geometry}
\usepackage{multicol}
\newtheorem{theorem}{Theorem}[section]
\newtheorem{lemma}{Lemma}[section]
\newtheorem{definition}{Definition}[section]
\usepackage{amsmath}

\title
{Large deviations for processes in random environments with jumps}

\author%[Ivan Matic]
{Ivan Matic\\
Mathematics Department, Duke University \\
Durham, NC 27708, USA\\
matic@math.duke.edu}

\begin{document}

\date{}

\maketitle
%\vspace{-1cm}

\begin{abstract} A deterministic walk in a random environment can be understood as a general random process with finite-range dependence  that starts repeating a loop once it reaches a site it has visited before. Such process lacks the Markov property. 
We study the exponential decay of the probabilities that the walk  will reach sites located far away from the origin. We also study a similar problem for the  continuous analogue: the  process that is a solution to an ODE with random coefficients. In this second model the environment also has ``teleports'' which are the regions from where the process can make discontinuous jumps. 
\end{abstract}


\vspace{0.5cm}


\vspace{0.5cm}

%\noindent Submitted to EJP on March 18, 2011, final version accepted November 14, 2011.
%\pagebreak

\section{Introduction}\label{in}
A deterministic process in a random environment is a solution to a differential equation of the form $\frac{dX_t}{dt}=b(X_t,\omega)$, $X_0=0$, where 
$b:\mathbb R^d\times \Omega\to\mathbb R^d$ is a vector field defined on some probability space $(\Omega, \mathbb P)$. We are interested in large deviation properties of the solution $X_t$.  This process is called deterministic because once the environment is sampled from the probability space $\Omega$, the entire process depends only on the initial position $X_0$. 

We will consider a modified version of the previous process, in which we assume certain randomness in time. The environment is placed in a Poisson point process, and small neighborhoods of the Poisson points serve as teleports. Once the particle spends sufficient time within a teleport, it is subject to a jump. However, whether a jump occurs or not depends on an additional time-dependent random sequence. 

One motivation for studying such processes is a simplified model for the spread of  viruses. Assume that the primary mode of infection spread is a local contact of members of a group. However, the mutation of other viruses into the virulent strain could appear in isolated areas where they are not introduced through direct contact. We model this phenomena by random jumps.

The model we study is also related to killed Brownian motions in Poisson potentials. The killed Brownian motion evolves according to a stochastic ODE, but once it spends too much time in a neighborhood of a Poisson point, it dies. Such processes are discussed in  \cite{sznitmanbook} where certain large deviation results are obtained. Our process continues to live, but reappears at another location. 

Under assumptions on $b$ and the environment, we will establish the following result about the moment generating functions of $X_t$.

\begin{theorem}\label{det_process}
Under the assumptions from Section \ref{section_dpre} there exists a convex 
function $\Lambda: \mathbb R^d\rightarrow\overline{\mathbb R}$ such that $$\lim_{t\to+\infty}\frac1t\log\mathbb E\left[e^{\lambda\cdot X_t}\right]=\Lambda(\lambda).$$
\end{theorem}

Here $\mathbb E$ denotes the expected value with respect to $(\Omega,\mathbb P)$. This is known as an {\em annealed} or {\em averaged} process.

To get an idea on how to approach the problem we first consider a discrete analogue -- a deterministic walk in a random environment. 
We consider a random sequence $X_n$ defined recursively as $X_{n+1}-X_n=b(X_n,\omega)$, for a suitable function $b:\mathbb Z^d\times \Omega\rightarrow \mathbb Z^d$.  
In this case we may assume that  $b(x,\omega)$ is of the form 
$b(x,\omega)=\eta_x$, where $\eta$ is a random field on $\mathbb Z^d$. The field $\eta$ itself could be understood as the random environment in which the walk occurs. When the particular realization of the environment is fixed, the walk becomes deterministic. 
Our aim is to establish the following theorem.

\begin{theorem}
Let $(\eta_z)_{z\in\mathbb Z^d}$ be a stationary $\mathbb Z^d$-valued random field that satisfies the assumptions (i)--(iii) (we refer to the assumptions from Section \ref{section_dwre}). 
Assume that the sequence $(X_n)_{n=1}^{\infty}$ of random variables is defined recursively as $X_0=0$, $X_{n+1}=X_n+\eta_{X_n}$. Then there exists a convex function $\Lambda: \mathbb R^d\rightarrow\overline{\mathbb R}$ such that $$\lim_{n\rightarrow\infty}\frac1n\log\mathbb E\left[e^{\lambda\cdot X_n}\right]=\Lambda(\lambda).$$
\end{theorem}

In order to make connections of this model with the one of the random walk in a random environment, we start from the representation of $X_n$  as $X_{n+1}-X_n=b(X_n,\omega,\pi)$. Here  $(\Omega, \mathbb P)$ and $(\Pi, P)$ are probability spaces and $b:\mathbb Z^d\times\Omega\times \Pi\rightarrow \mathbb Z^d$ 
is a random variable. In this model, $(\Omega, \mathbb P)$ is the environment, and for each fixed $\omega\in\Omega$, the walk $X_n$ could be understood as a random walk on probability space $(\Pi, P)$. Recent works  \cite{varadhan} and \cite{rassoul-agha} have established the quenched and annealed large deviations for the random walk in a random environment under certain non-degeneracy assumptions on $(\Pi, P)$. The articles \cite{rosenbluth} and \cite{yilmaz} found variational formulas for the rate functions and some connections between the two rate functions. The model we are studying is related to the annealed (averaged) case studied in the mentioned papers. Our probability space $(\Pi, P)$ is very degenerate.

The main idea in our study is the adaptation of the subadditive argument. However, we need some work to find a quantity that happens to be subadditive. This quantity is found in  Section \ref{hitting_times} where we prove that it has the same exponential rate of decay as our large deviation probabilities. After that we are ready for the proofs of the main theorems which are presented in Section \ref{section_ld}. Section \ref{section_lln} deals with consequences of large deviations such as the law of large numbers and the existence of a loop. The law of large numbers with limiting velocity 0 for the deterministic walk in random environment is not very surprising. Intuitively, the expected time until the occurrence of the loop is finite which means that the walk can't go too far.

We will restrict ourselves to environments with finite-range dependence. The special cases are iid environments where many of the presented arguments can be simplified. We are not assuming that the walk is nearest--neighbor. We do assume however, certain {\em ellipticity} conditions on the environment: There is a {\em nice} set of vectors such that using only these vectors the walk can move in any direction. Our ellipticity condition is the requirement that the probability that the walk at any position of the environment takes any particular vector from the  nice set is uniformly positive. The existence of the nice set guarantees that there are no forbidden directions for the walk. 

In the next section we state the assumptions we impose to the discrete model and state 
Theorem \ref{main_theorem} which is the main focus of our study. In order to prove it we will need some deterministic lemmas that establish uniform bounds on probabilities that the walk can go in any direction. These statements are stated and proved in Section \ref{existence_of_loop}. Subadditive quantities are found in Section \ref{hitting_times} and in Section \ref{section_ld} we use them to prove Theorem \ref{mainmain_theorem}. The law of large numbers and the existence of a loop are proved in Section \ref{section_lln}.
In the last section we discuss the generalization of this approach to the continuous setting.

Many of the arguments here involve the consideration of the position of the particle at the time of the first self-intersection. Other interesting results regarding self-intersection times of random walks can be found  in \cite{aldous} and \cite{asselah}.
The hitting times and back-tracking times were a useful tool in establishing the law of large numbers for ballistic random diffusions (see \cite{shen}).

\section{Definitions and Assumptions}
Throughout the paper we will use $\mathbb R_+$ to denote the set of positive real numbers and 
$\bar{\mathbb R}=\mathbb R\cup \{-\infty,+\infty\}$. For $v\in\mathbb R^d$ we will use $|v|$ to denote its 
norm. 
The indicator function for a measurable set 
$B$ will be denoted by  $1(B)$.  The ball with center $x$ and radius $r$ will be 
denoted by $B_r(x)$, or $B_r$ for brevity if the center is clear from the context. 
For $x\in\mathbb R$ we will denote by $\lfloor x\rfloor$ (resp. $\lceil x\rceil$) the 
largest integer not greater than $x$ (resp. the smallest integer not smaller than $x$). 
The complement of a set $Q$ will be denoted by $Q^C$.

\label{section_dwre}
\begin{definition}
The set of vectors $\{u_1, \dots, u_m\}\subseteq\mathbb Z^d$ is called {\em nice} if for every $l\in\mathbb R^d\setminus\{0\}$ 
there exists $i\in\{1,2,\dots,m\}$ such that $l\cdot u_i>0$. 
\end{definition}
Let $(\eta_z)_{z\in\mathbb Z^d}$ be a stationary $\mathbb Z^d$-valued random field that satisfies  the following conditions: 
\begin{enumerate}
\item[(i)] There exists a positive real number $L$ such that $|\eta_z|\leq L$ for all $z\in\mathbb Z^d$. 
\item[(ii)] There exists a real number $M$ such that $\eta_{z}$ is independent 
of the environment outside of the ball with center $z$ and radius $M$. 
\item[(iii)] There exist a nice set of vectors $\{u_1, \dots, u_m\}\in\mathbb Z^d$ and  a constant $c>0$ such that $\mathbb P(\eta_z=u_i|\mathcal F_z)>c$ for all $i\in\{1,2,\dots, m\}$, where $\mathcal F_z$ is 
a sigma-algebra generated by all $\eta_{w}$ for $w\in\mathbb Z^d$ such that $0<|w-z|\leq M$.
\end{enumerate} 
The last assumption implies the existence of a loop in any half-space with a positive 
probability (see Theorem \ref{theorem_1}). It also implies that there exists a constant 
$c>0$ such that $\mathbb P(\eta_z\cdot l>0|\mathcal F_z)>c$ for every $z$ and every $l\in
\mathbb R^d\setminus\{0\}$.

A special case is the iid environment when in condition (ii) we require $M<1$.  The condition (iii) is then replaced by $\mathbb P(\eta_z=u_i)>c$ for all $i$. 

The random variable $X_n$ is defined recursively as $X_0=0$ and $$X_{n+1}=X_n+\eta_{X_n}\;\;\mbox{ for }\;\; n\geq 1.$$  

We will use the following equivalent interpretation of $X_n$. The process $X_n$ behaves like a random walk until the first self intersection. The increments of the random walk are sampled at every step according to the law of the random field $\eta_z$. After the self intersection occurs, the walk becomes deterministic and repeats the loop. 

Here is the precise definition of the walk $X_n$ that we will use. Let $\xi_n$ be a $\mathbb Z^d$-valued random sequence defined recursively:  $\xi_0=0$, and 
\begin{eqnarray}\label{definition_xi} \xi_n&=&\eta_{\xi_0+\cdots + \xi_{n-1}} \;\;
\mbox{ for }\;\;n\geq 1.
\end{eqnarray}
 Let us define $Y_n=\sum_{k=1}^n \xi_k$ and 
$\tau=\inf\{n: Y_n\in\{Y_1,\dots, Y_{n-1}\}\}$. 
If $\tau<\infty$, let $\theta< \tau$ be the smallest integer such that $Y_{\tau}=Y_{\theta}$. 
Denote by $\theta'$ the remainder when $n-\theta$ is divided by $\tau-\theta$. Formally,
 $\theta'=(n-\theta)-(\tau-\theta)\lfloor (n-\theta)/(\tau-\theta)\rfloor$. 
This remainder will belong to the set $\{0,1,\dots, \tau-\theta-1\}$.

We  define the walk $X_n$ using the formula:
$$X_n=\left\{
\begin{array}{ll}Y_n, &n\leq \tau \mbox{ or }\tau=\infty,\\
Y_{\theta+\theta'}, & \mbox{otherwise}.\end{array}
\right.
$$

Let $l\in\mathbb R^d$ be a unit vector. Define $T_m^l=\inf\{n: X_n\cdot l\geq m\}$.
For $x\in\mathbb R^d$ let us  
denote by $H_x^l$ the half-space through $x$ determined by the vector $l$ as $$H_x^l=\mathbb Z^d\cap\left\{x+v:v\in\mathbb R^d, l\cdot v\geq 0\right\}.$$
For $s\in\mathbb R$ denote by $Z_s^l$ the hyperplane
$$Z_s^l=\{x\in\mathbb R^d: x\cdot l=s\}.$$

Our goal is to prove the large deviations for $X_n$ (see Theorem \ref{mainmain_theorem}). 
We will be able to use the Ellis-G\" artner theorem to get some further bounds once we establish the following result:

\begin{theorem}\label{main_theorem} Let $X_n$ be the random walk defined as above. Assume that the random environment satisfies the conditions (i)--(iii). 
For each unit vector $l\in\mathbb R^d$ there exists a concave function $\phi^l:\mathbb R_+\rightarrow\bar{\mathbb R}$ such that for all $k\in\mathbb R_+$:
\begin{eqnarray}\label{main_limit}
\lim_{n\rightarrow\infty} \frac1n\log\mathbb P(X_n\cdot l\geq nk)&=&\phi^l(k).
\end{eqnarray}
\end{theorem}

\noindent{\em Remark.} Notice that $\phi^l(k)=\phi^{tl}(tk)$ for all $t\in\mathbb R_+$. Therefore $\phi^l(k)=\Phi\left(\frac1kl\right)$ for a suitable function $\Phi:\mathbb R^d\rightarrow\bar{\mathbb R}$. Notice that the concavity or even continuity of 
$\Phi$ does not follow from the properties of $\phi^l$. One expects that $\Phi$ is continuous, but we are not able to prove that, and this is one of the major unsolved problems. In the end we will discuss a bit about difficulties we are facing in trying to address the continuity of $\Phi$.

\section{Existence of a Loop}
\label{existence_of_loop}

In this section we prove that the previously defined random walk will have a loop in each half-space with a positive probability. This fact will be a consequence of the following elementary lemma. The lemma states that there exists a loop consisting entirely of vectors from a nice set. 

\begin{lemma}\label{lemma_1}
Let $\{u_1$, $\dots$, $u_m\}$ be a nice set of non-zero vectors. There exist non-negative integers $q_1$, $q_2$, $\dots$, $q_m$ not all equal to $0$ such that $q_1u_1+\cdots + q_mu_m=0$.
\end{lemma}
\noindent{\bf Proof.}
We will prove the statement using the induction on the dimension $d$. The statement is easy to prove for $d=1$ and $d=2$.
We may assume that $\{u_1,\dots, u_m\}$ is a minimal nice set, i.e. there is no proper nice subset of $\{u_1, \dots, u_m\}$. If not, take the proper nice subset and repeat the argument. Let us fix the vector $u_m$, and let $$v_i=u_i-\frac{u_i\cdot u_m}{|u_m|^2}u_m\; \mbox{ for }\; i=1, \dots, m-1.$$ All vectors $v_1$, $\dots$, $v_{m-1}$ have rational coordinates. Let $r$ be the common denominator of those fractions and consider the lattice $D$ of size $1/r$ in the vector space determined by the span $W$ of $v_1, \dots, v_{m-1}$. Let us prove that the set $\{v_1, \dots, v_{m-1}\}$ is nice in $W$. Let $\tilde l\in W$ be a vector with real coordinates. There exists $i\in\{1,2,\dots, m\}$ such that $u_i\cdot \tilde l>0$. 
Since $\tilde l\in W$ we immediately have that $u_m\cdot \tilde l=0$ and $u_i\cdot \tilde l=
v_i\cdot \tilde l$ hence $v_i\cdot \tilde l>0$. This implies that $\{v_1, \dots, v_{m-1}\}$ is a nice set of vectors in $W$. 
According to the induction hypothesis there are non-negative integers $q_1', \dots, q_{m-1}'$ such that $q_1'v_1+\cdots + q_{m-1}'v_{m-1}=0$.  We now have $$|u_m|^2(q_1'u_1+\cdots + q_{m-1}'u_{m-1})=\left(q_1'u_1\cdot u_m+\cdots+q_{m-1}'u_{m-1}\cdot u_{m}\right)u_m.$$
Let us now prove that $q_1'u_1\cdot u_m+\cdots + q_{m-1}'u_{m-1}\cdot u_m\leq 0$. Assume the 
contrary, that this number were greater than $0$. Since $\{u_1, \dots, u_{m-1}\}$ is 
not a nice set (due to our minimality assumption for $\{u_1,\dots, u_m\}$) there exists a vector $l\in\mathbb R^d$ such that $l\cdot u_m>0$ but $l\cdot u_k\leq 0$ for each $k\in\{1,2,\dots, m-1\}$. 
This gives that 
\begin{eqnarray*}0&\geq&l\cdot (q_1'u_1+\cdots + q_{m-1}'u_{m-1})|u_m|^2
=\left(q_1'u_1\cdot u_m+\cdots + q_{m-1}'u_{m-1}\cdot u_m\right) u_m\cdot l>0,
\end{eqnarray*}
a contradiction. 
Therefore $q_1'u_1\cdot u_m+\cdots+ q_{m-1}'u_{m-1}\cdot u_m\leq 0$. We can now choose $q_i=|u_m|^2q_i'$, $i=1,2,\dots, m-1$, and $q_m=-(q_1'u_1\cdot u_m+\cdots+ q_{m-1}'u_{m-1}\cdot u_m)$ to obtain $q_1u_1+\cdots + q_mu_m=0$. This completes the proof.
\hfill $\Box$

\vspace{0.5cm}

The paper \cite{rassoulagha_seppalainen_yilmaz} contains a proof that each point of the lattice can 
be accessed by making only moves governed by the vectors $u_1$, $\dots$, $u_m$. 

The following theorem says that in each half-space $H^l_x$  a walk $X_n$ starting from $x$ can have a loop in $H^l_x$ with a probability that is strictly greater than $0$. Because of the stationarity it suffices to prove this for half-spaces through the origin.

\begin{theorem}\label{theorem_1}

Let $U=\{u_1,\dots, u_m\}$ be a nice set of non-zero vectors. 
There exists a constant $m_0\in\mathbb N$ such that:
For each unit vector $l\in\mathbb R^d\setminus\{0\}$ there  exist an integer $s\leq m_0$ 
and a sequence $x_0=0$, $x_1$, $x_2$, $\dots$, $x_s$ $\in U$ with the following
properties:
\begin{enumerate}
\item[(i)] $\sum_{i=0}^p x_i\cdot l\geq 0$ for all $0\leq p\leq s$,
\item[(ii)]
$\sum_{i=0}^px_i\neq \sum_{i=0}^q x_i$ for all $0\leq p<q\leq s-1$, and 
\item[(iii)] $\sum_{i=0}^s x_i=0$.
\end{enumerate}
\end{theorem}
\noindent{\bf Proof.} 
Using Lemma \ref{lemma_1} we have that there exists a sequence 
$z_1, \dots, z_s\in\{u_1, \dots, u_m\}$ such that $z_1+\cdots + z_s=0$. Let us point out that 
$z_i$ need not be distinct. 
Assume that this sequence is one of the sequences of minimal length. 
Let us choose the index $\pi$ such that that $(z_1+\cdots + z_{\pi})\cdot l$ is 
minimal (of all $j=1,2,\dots, s$). Then 
$x_1=z_{\pi+1}$, $x_2=y_{\pi+2}$, $\dots$, $x_s=y_{\pi+s}$ 
(indices are modulo $s$) satisfy $x_1+x_2+\cdots + x_s=0$, and 
$$(x_1+x_2+\cdots + x_{j})\cdot l\geq 0\;\;\mbox{ for all }j=1,2,\dots, s.$$ 
Let us prove the last inequality. Assume the contrary, that 
$(x_1+\cdots + x_{j})\cdot l<0$ for some $j$. Assume first that $\pi+j\leq s$. Then 
$(z_1+\cdots + z_{\pi+j})\cdot l=(z_1+\cdots + z_{\pi})\cdot l + (x_1+\cdots + x_{j})
\cdot l<(z_1+\cdots + z_{\pi})\cdot l$ which contradicts the choice of $\pi$. If 
$\pi+j>s$ and $\pi+j\equiv t$ (mod $s$) then 
$0=l\cdot (z_1+\cdots + z_{\pi} + x_1+\cdots + x_{s-\pi})<
l\cdot (z_1+\cdots + z_{t}+x_1+\cdots + x_{s-\pi})=l\cdot (x_1+\cdots + x_j)<0$, 
a contradiction. This proves (i). The condition (ii) follows from the requirement 
that the sequence $z_1, \dots, z_s$ is the shortest. \hfill $\Box$

\vspace{0.5cm}

Let us recall the definition of $\xi_i$ from (\ref{definition_xi}). The following result 
establishes a useful property of the sequence $x_0$, $\dots$, $x_s$ that we constructed in the  
previous theorem.
\begin{theorem}\label{theorem_11}
Let $y_i=x_0+\cdots + x_{i-1}$ for $i\in\{1,2,\dots, s\}$ and let us 
denote by $\mathcal F_{y_1, \dots, y_s}$ the $\sigma$-algebra generated by all 
random variables $\eta_z$ for $z\in\mathbb Z^d\setminus  \{y_1, \dots, y_s\}$ such 
that $ \min_{1\leq i\leq s}|z-y_i|\leq M$. Then:
$$\mathbb P(\xi_{1}=x_1, \xi_{2}=x_2, \dots, 
 \xi_{s}=x_s|\mathcal F_{y_1, \dots, y_s})\geq c_1. $$
\end{theorem}
\noindent{\bf Proof.}
Denote by $\mathcal G_{y_s}^{y_1, \dots, y_{s-1}}$ the sigma algebra generated by all random variables $\eta_z$ for $z\in \mathbb Z^d\setminus \{y_s\}$ such that $\min_{1\leq i\leq s}|z-y_i|\leq M$. 
\begin{eqnarray*}&&
\mathbb P(\xi_1=x_1, \xi_2=x_2, \dots, 
 \xi_s=x_s|\mathcal F_{y_1, \dots, y_s})\\&=&\mathbb E\left(
\left. \mathbb E\left(
1( \xi_1=x_1) \cdots 1( 
 \xi_s=x_s)|\mathcal G_{y_s}^{y_1, \dots, y_{s-1}}\right)
\right|\mathcal F_{y_1, \dots, y_s}
 \right)\\&=&
 \mathbb E\left(1( \xi_1=x_1) \cdots 1( 
 \xi_{s-1}=x_{s-1})\cdot
\left. \mathbb E\left(1(\xi_{s}=x_s)
|\mathcal G_{y_s}^{y_1, \dots, y_{s-1}}\right)
\right|\mathcal F_{y_1, \dots, y_s}
 \right)\\&>&c\cdot \mathbb P(\xi_1=x_1, \dots, \xi_{s-1}=x_{s-1}|\mathcal F_{y_1,\dots, y_s})\\&=&c\cdot \mathbb P( \mathbb P(\xi_1=x_1, \dots, \xi_{s-1}=x_{s-1}|\mathcal F_{y_1,\dots, y_{s-1}})|\mathcal F_{y_1,\dots, y_{s}}).
\end{eqnarray*} 
Now we can continue by induction to obtain that 
$$\mathbb P(\xi_1=x_1, \dots, \xi_s=x_s|\mathcal F_{y_1, \dots, y_s})>c^s.$$
It remains to take $c_1=c^s$. \hfill $\Box$

\vspace{0.5cm}

Using the previous two lemmas we can establish the equality analogous to the one from Theorem \ref{main_theorem} in which $k=0$. 

\begin{theorem}\label{case_k=0} For each vector $l\in\mathbb R^d$ and $X_n$ defined as before the following equality holds:
$$\lim_{n\rightarrow\infty}\frac1n\log\mathbb P(X_n\cdot l\geq 0) =0.$$
\end{theorem}
\noindent{\bf Proof.} The inequality $\mathbb P(X_n\cdot l\geq 0)\leq 1$ implies that $\limsup\frac1n\log\mathbb P(X_n\cdot l\geq 0 )\leq 0$. For the other inequality we use Theorem \ref{theorem_1}. Let $x_1, \dots, x_s$ be the sequence whose existence is claimed by that theorem. Let $y_i=x_0+\cdots + x_{i-1}$ as before. Notice that 
\begin{eqnarray*}\mathbb P(X_n\cdot l\geq 0 ) &\geq &
\mathbb P(\xi_1=x_1, \xi_2=x_2, \dots, \xi_s=x_s)\\&=&
\mathbb E\left(\mathbb P\left(\left. \xi_1=x_1, \xi_2=x_2, \dots, \xi_s=x_s\right|\mathcal F_{y_1, \dots, y_s}\right)\right)\\&\geq & c_1.
\end{eqnarray*}
Therefore $\liminf\frac1n\log\mathbb P(X_n\cdot l\geq 0)\geq \liminf \frac1n\log c_1=0$.
\hfill $\Box$

\vspace{0.5cm}

We will also need the following deterministic lemma. 
\begin{lemma}\label{uniformity_nice}
Assume that $\{u_1, \dots, u_m\}\subseteq \mathbb Z^d$ is a nice set of vectors. Let $\rho:\mathbb R^d\rightarrow\mathbb R_+\cup\{0\}$ be the function defined as $\rho(l)=\max_i\{u_i\cdot l\}$. Then $$\inf_{l:|l|=1} \rho(l)>0.$$
\end{lemma}
\noindent{\bf Proof.} First notice that $\rho(l)>0$ for each $l\in\mathbb R^d\setminus\{0\}$. Otherwise the set $\{u_1, \dots, u_m\}$ would not be nice. Notice also that $\rho$ is a continuous function (because it is a maximum of $m$ continuous functions) and the unit sphere is a compact set. Thus the infimum of $\rho$ over the unit sphere must be attained at some point, and we have just proved that value of $\rho$ at any single point is not $0$. \hfill $\Box$

\section{Hitting Times of  Hyperplanes}\label{hitting_times}

The main idea for the proof of Theorem \ref{main_theorem} is to establish the asymptotic equivalence of the sequence $\frac1n\log\mathbb P(X_n\cdot l\geq nk)$ and a sequence to which we can apply the deterministic superadditive lemmas. First we will prove that the previous sequence behaves  like $\frac1n \log\mathbb P(T_{nk}^l\leq n)$. Then we will see that the asymptotic behavior of the latter sequence satisfies $$\frac1n \log\mathbb P(T_{nk}^l\leq n)\sim
\frac1n\log\mathbb P(T_{nk}^l\leq n, T_{nk}^l\leq D_1^l)$$ where $D_1^l$ is 
the first time the walk returns over the hyperplane $Z_0^l$. We  define $D_1^l$ precisely in 
the following way. First, let $B_1^l=\inf\{n>0: X_n\cdot l< 0\}$. Then
$D_1^l=\inf\{n>B_1^l: X_n\cdot l\geq 0\}$. We allow $D_1^l$ to be infinite.
The event $\{T_{nk}^l\leq n, T_{nk}^l\leq D_1^l\}$ describes 
those walks that reach the hyperplane $Z_{nk}^l$ without  backtracking 
over the hyperplane $Z_0^l$.

We will be able to prove the existence of the limit of the last 
sequence using a modification of the standard subadditive result that states that $\lim \frac{a_n}n=\inf\frac{a_n}n$ if $a_{n+m}\leq a_n+a_m$ for all $m,n\in\mathbb N$. 

From now on we fix the unit vector $l\in\mathbb R^d$ and  omit the 
superscript $l$ in the variables. Also, some of the new variables that will be 
defined would need to have a superscript $l$ but we will omit it as well. 

Our first result in carrying out the formerly described plan  is the following lemma: 

\begin{lemma}\label{first_lemma}
The following inequality holds: 
\begin{eqnarray*}
\limsup_{n\rightarrow \infty}\frac1n \log \mathbb P(X_n\cdot l\geq nk)&\leq&
\limsup_{n\rightarrow\infty} \frac1n \log\mathbb P(T_{nk}\leq n)
\end{eqnarray*}
In addition, for each $\varepsilon>0$ we have:
\begin{eqnarray*}
\liminf_{n\rightarrow\infty} \frac1n \log\mathbb P(T_{n(k+\varepsilon)}\leq n)&\leq&
\liminf_{n\rightarrow \infty}\frac1n \log \mathbb P(X_n\cdot l\geq nk).
\end{eqnarray*}
\end{lemma} 

\noindent{\bf Proof.}
Clearly, $\{X_n\cdot l \geq kn\}\subseteq \{T_{kn}\leq n\}$. Therefore $\mathbb P(X_n\cdot l\geq kn)\leq \mathbb P(T_{kn}\leq n)$. This establishes the first inequality. Let $x_0, \dots, x_s$ be the sequence whose existence follows from Theorem \ref{theorem_1}.
We now have
\begin{eqnarray} \nonumber
\mathbb P(X_n\cdot l\geq kn)&=& \mathbb P(X_n\cdot l\geq kn, T_{kn}\leq n)\\
\nonumber &\geq&
\mathbb P\left( X_n\cdot l \geq kn, T_{kn}\leq n-s, X_{T_{kn}+1}-X_{T_{kn}}=x_1, 
\dots,
%\right.\\
%\nonumber
%&&\left.
 X_{T_{kn}+s}-X_{T_{kn}+s-1}=x_s\right)
\\ \label{forget_xncdotl}
&=& % \nonumber
\mathbb P\left( T_{kn}\leq n-s, X_{T_{kn}+1}-X_{T_{kn}}=x_1, 
\dots,
%\right.\\&&\left. 
X_{T_{kn}+s}-X_{T_{kn}+s-1}=x_s\right)\\ \nonumber
&=&\mathbb E\left[1(T_{nk}\leq n-s)\mathbb E\left[1( \xi_{T_{kn}}=x_1) \cdots
%\right.\right.\\ \nonumber&&\left.\left.
\left.  1(\xi_{T_{kn}+s-1}=x_s)\right|\mathcal F_{T_{kn}}\right]\right]\\ \nonumber
&\geq &c_1\cdot \mathbb P\left(T_{kn}\leq n-s
\right).
\end{eqnarray}
Here $\mathcal F_{T_{kn}}$ denotes the $\sigma$-algebra defined by $\eta_z$ for $z\in\mathbb Z^d$ such that $|z-X_i|\leq M$ for $i=1,2,\dots, T_{kn}$. The equality in (\ref{forget_xncdotl})  holds because if $T_{kn}\leq n-s$ and $X_{T_{kn}+1}-X_{T_{kn}}=x_1$, $\dots$, 
$X_{T_{kn}+s}-X_{T_{kn}+s-1}=x_s$, then the walk will enter in a loop. This loop will be in the half-space $H_{kn}$ which would guarantee that $X_n\cdot l>kn$. 

For each $\varepsilon>0$, if $n>s+sk/\varepsilon $ we have  
$\{T_{(k+\varepsilon)(n-s)}\leq n-s\}\subseteq \{T_{kn}\leq n-s\}$, and consequently the first set has the smaller probability than the second.
This completes the proof of the lemma. \hfill $\Box$

\vspace{0.5cm}
 
\noindent{\em Remark.} In  the same way we could obtain the analogous inequalities with the walk $X_n$ replaced by
$X_{n\wedge \tau}$. 

 For each integer $i\geq 1$ denote by $D_i$ the time of the 
$i$th jump over the hyperplane $Z_0$ in the direction of the vector $l$. 
Define $D_0=0$, $B_i=\inf\{n>D_{i-1}:X_n\cdot l<0\}$ and 
$D_i=\inf\{n>B_i: X_n\cdot l\geq 0\}$. We allow for $D_i$ to be $\infty$. 

\begin{lemma}\label{second_lemma}
Let $k$ and $k'$ be two real numbers such that $0<k'<k$. Then the following two inequalities hold: 
\begin{eqnarray*}\limsup \frac1n\log\mathbb P(T_{nk}\leq n) &\leq& \limsup\frac1n\log\mathbb P(T_{nk'}\leq n, T_{nk'}\leq D_1)\\
\liminf \frac1n\log\mathbb P\left(T_{nk}\leq n, T_{nk}\leq D_1\right)&\leq &\liminf
\frac1n\log\mathbb P\left(T_{nk}\leq n\right).\end{eqnarray*}
\end{lemma}  
\noindent{\bf Proof.} 
The second inequality is easy and it follows from  
 $$\{T_{nk}\leq n\}\supseteq \{T_{nk}\leq n, T_{nk}\leq D_1\}.$$
Our goal is to prove the first inequality. We start by noticing that 
\begin{eqnarray*}
\mathbb P\left(T_{u}\leq n\right)&=&
\sum_{i=1}^n \mathbb P\left(T_{u}\leq n, 
D_i< T_{u}\leq D_{i+1}\right).
\end{eqnarray*}
We will prove that each term from the right-hand side of the previous equality is bounded above by $L^{d}n^{d-1}\cdot\mathbb P(T_{u-L}\leq n, T_{u-L}\leq D_1)$. 
Let $Z_0'$ be the set of all points $z\in\mathbb Z^d$ such that $z\cdot l>0$ and the distance between $z$ and $Z_0$ is at most $L$. We have:
\begin{eqnarray*}
\mathbb P\left(T_{u}\leq n, 
D_i< T_{u}\leq D_{i+1}\right)
&=&\sum_{z\in  Z_0', |z|\leq nL} 
\mathbb E\left[1(T_{u}\leq n)\cdot 1(X_{D_i}=z)\cdot1(T_{u}> D_i)\right.\\&&\left.\cdot 1(T_{u}\leq D_{i+1})
\right]\\
&=&\sum_{z\in  Z_0', |z|\leq nL} \mathbb E\left[
\mathbb E\left[1(T_{u}\leq n)\cdot 1(X_{D_i}=z)\cdot1(T_{u}> D_i)\right.\right.\\&&\left.\left.\cdot\left. 1(T_{u}\leq D_{i+1})
\right|\mathcal F_{D_i,T_u}\right]\right].
\end{eqnarray*}
Here  $\mathcal F_{D_i,T_u}$ denotes the $\sigma$-algebra generated by the 
random environment that is contained in the $M$-neighborhood of the walk 
from $D_i$ to $T_u$. When conditioning on this $\sigma$-algebra we 
essentially understand  our environment in the following way: 
It consists of two walks: One deterministic that goes from $z$ to $Z_{u}$ 
without crossing the hyperplane $Z_0$,  and another walk  that starts at $0$ 
and ends in $z$ by  making exactly $i$ crossings over $Z_0$, not intersecting 
the other deterministic walk, and not crossing over $Z_u$. Therefore: 
\begin{eqnarray*}
&&\sum_{z\in  Z_0', |z|\leq nL} \mathbb E\left[
\mathbb E\left[1(T_{u}\leq n)\cdot  1(X_{D_i}=z)\cdot1(T_{u}> D_i)\cdot 1(T_{u}\leq D_{i+1})\left.
\right|\mathcal F_{D_i,T_u}\right]\right]\\
&\leq &\sum_{z\in  Z_0', |z|\leq nL} \mathbb E\left[
1(\tilde T_{u-L}\leq n-D_i) \cdot 1(\tilde T_{u-L}\leq\tilde D_1)
 \cdot 
\mathbb E\left[1(X_{D_i}=z)\cdot1(T_{u}> D_i)|\mathcal F_{D_i,T_u}\right]\right]\\
&\leq & \sum_{z\in  Z_0, |z|\leq nL}
\mathbb E\left[1(\tilde T_{u-L}\leq n)\cdot 1(\tilde T_{u-L}\leq \tilde D_1)
\right],
\end{eqnarray*}
where $\tilde T_u$ is defined in analogous way as $T_u$ to correspond to the new 
walk $\tilde X_j=X_{D_i+j}$. 
In the last equation $\tilde D_1$ is defined as the first time of crossing over 
the hyperplane $\Gamma$ parallel to $Z_0$ that is shifted by the vector 
$-Ll$. Let us now prove that 
$\mathbb P(\tilde T_{u-L}\leq n,\tilde T_{u-L}\leq \tilde D_1)\leq L^dn^{d-1}\cdot 
\mathbb P(T_{u-2L}\leq n, T_{u-2L}\leq D_1)$. Denote by $J$ the 
time when the walk comes closest to the hyperplane $\Gamma$. 
For the precise definition of $J$ we first introduce $\mu=\inf\{
\mbox{dist }(X_n,\Gamma):D_i<n\leq D_{i+1}\}$. Here 
$\mbox{dist }(X,\pi)$ denotes the distance between the point $X$ and the hyperplane $\pi$. 
We now define $J=\inf\{n: \mbox{ dist }(X_n,\Gamma)=\mu\}$.
The number of possibilities for $X_j$ is at most $L^dn^{d-1}$ and similarly as 
above, conditioning on the $\sigma$-algebra between $J$ and $T_{u-L}$ we get  
$\mathbb P( T_{u-L}\leq n, T_{u-L}\leq \tilde D_1)\leq L^dn^{d-1}\cdot 
\mathbb P(T_{u-2L}\leq n, T_{u-2L}\leq D_1)$. 
We now have  $$\mathbb P(T_u\leq n, D_i\leq T_u\leq D_{i+1})\leq L^{2d}n^{2d-2}\mathbb P(T_{u-2L}\leq n, T_{u-2L}\leq D_1).$$ This implies that $\mathbb P(T_u\leq n)\leq L^{2d}n^{2d-1}\mathbb P(T_{u-2L}\leq n, T_{u-2L}\leq D_1)$.
Therefore
$$\limsup \frac1n\log\mathbb P(T_{nk}\leq n) \leq \limsup\frac1n\log\mathbb P(T_{nk-2L}\leq n, T_{nk-2L}\leq D_1).$$
For each $k'<k$ we have
$$\limsup \frac1n\log\mathbb P(T_{nk}\leq n)\leq \limsup\frac1n\log\mathbb P(
T_{nk'}\leq n, T_{nk'}\leq D_1).$$
This completes the proof of the lemma. \hfill $\Box$

\section{Large Deviations Estimates}\label{section_ld}
Now we are able to prove the main theorem:

\noindent
{\bf Proof of Theorem \ref{main_theorem}.} 
We will prove that for each unit vector $l$ there exists a concave function $\psi:\mathbb R_+\rightarrow\bar{\mathbb R}$ such that 
\begin{eqnarray}\label{limit_psi}\lim_{n\rightarrow\infty}\frac1n\log\mathbb P(T_{nk}\leq n) &=& \psi(k).\end{eqnarray}
Because of  Lemma \ref{second_lemma} it suffices to prove 
that there exists a concave function $\gamma:\mathbb R_+\rightarrow\bar{\mathbb R}$ such that 
\begin{eqnarray}
\label{limit_gamma}
\lim_{n\rightarrow\infty}\frac1n\log\mathbb P\left(T_{nk}\leq n, T_{nk}\leq D_1\right)&=&\gamma(k).\end{eqnarray}

Let $w\in\mathbb Z^d$ be a vector such that $w\cdot l>0$ and $\mathbb P(\eta_z=w|\mathcal F_z)\geq c>0$ for some constant $c$. Assume that $r$ is an integer such that 
the distance between the hyperplanes  $Z_{rw}$ and $Z_0$ is at least $M$. 
If $u,v,p,q$ are any four positive real numbers such that $q>r$, then the 
following inequality holds: 
\begin{eqnarray}
%\nonumber
\mathbb P\left(T_{u+v}\leq p+q, T_{u+v}\leq D_1\right)&\geq& c\cdot \mathbb P(T_u\leq p, 
T_u\leq D_1)
%\\&&
\label{sub_add}
\cdot \mathbb P(T_v\leq q-r, T_v\leq D_1).
\end{eqnarray}
If the environment were iid this could have been done by conditioning on $\mathcal F_u$. In our situation the idea is the same,  we just need some more work to compensate for the lack of independence.
\begin{eqnarray*}
\mathbb P\left(T_{u+v}\leq p+q, T_{u+v}\leq D_1\right)&\geq&\mathbb
P\left(T_{u+v}\leq p+q, T_{u+v}\leq D_1, T_u\leq p\right)\\
&\geq&\mathbb E\left[1(T_{u}\leq p)\cdot 1(T_u\leq D_1)\cdot 1(\xi_{T_u}=w)\cdot \cdots\right.\\
&&\left. \cdot 1(\xi_{T_u+r-1}=w)\cdot 1(T_{u+v}\leq p+q)\cdot 1(T_{u+v}\leq D_1)\right]
\\&=&\mathbb E\left[\mathbb E\left[1(T_{u}\leq p)\cdot 1(T_u\leq D_1)\cdot 1(\xi_{T_u}=w)\cdot \cdots  \cdot 1(\xi_{T_u+r-1}=w)\right.\right.\\
&&\left.\left.\left.\cdot 1(T_{u+v}\leq p+q)\cdot 1(T_{u+v}\leq D_1)\right|\mathcal F_{\left[X_{T_u},X_{T_{u+rw\cdot l}}\right]}\right]\right],
\end{eqnarray*}
where $\mathcal F_{[a,b]}$ is the $\sigma$-algebra determined by the 
environment outside of the strip $[Z_a,Z_b]$. Let us introduce the following 
notation: $\hat X_i=X_{T_u+r-1+i}$,  $\hat D_1$ the first time 
$\hat X_i$ jumps over $Z_{X_{T_u}+rw}$, and $\hat T_v=\inf\{i: (\hat X_i-\hat X_0)\cdot l\geq v\}$.
We now have that \begin{eqnarray*}
\mathbb P(T_{u+v}\leq p+q, T_{u+v}\leq D_1)&\geq & 
\mathbb E\left[\mathbb E\left[1(T_{u}\leq p)\cdot 1(T_u\leq D_1)\cdot 1(\xi_{T_u}=w)\cdot \cdots  \cdot 1(\xi_{T_u+r-1}=w)\right.\right.\\
&&\left.\left.\left.\cdot 1(\hat T_{v}\leq q-r)\cdot 1(\hat T_{v}\leq \hat D_1)\right|\mathcal F_{\left[X_{T_u},X_{T_{u+rw\cdot l}}\right]}\right]\right]\\
&=&
\mathbb E\left[1(T_{u}\leq p)\cdot 1(T_u\leq D_1)\cdot 1(\hat T_{v}\leq q-r)\cdot 1(\hat T_{v}\leq \hat D_1)\cdot \right.\\&&\left.
\mathbb E\left[\left. 1(\xi_{T_u}=w)\cdot \cdots  \cdot 1(\xi_{T_u+r-1}=w)\right|\mathcal F_{\left[X_{T_u},X_{T_{u+rw\cdot l}}\right]}\right]\right]\\
&\geq & c\cdot \mathbb P(T_u\leq p, T_u\leq D_1)\cdot \mathbb P(T_v\leq q-r, T_v\leq D_1).
\end{eqnarray*}
Let us explain why one can bound  $\mathbb E[1(\xi_{T_u}=w)\cdots 
1(\xi_{T_u+r-1}=w)|\mathcal F_{[X_{T_u}, X_{T_u+rw\cdot l}]}]$ below by $c$. 
The quantity in question is a random variable, and we are only considering that 
random variable on the set $\{T_u\leq D_1\}\cap \{\hat T_v\leq \hat D_1\}$. On that 
set the random walk doesn't visit the sites $X_{T_u}$, $\dots$, $X_{T_u+r-1}$ hence we can use our ellipticity assumption. 
This establishes (\ref{sub_add}). Let $\delta(m,x)=\log\mathbb P(T_{mx}\leq m, T_{mx}\leq D_1)$.
Applying the inequality (\ref{sub_add}) to the numbers $u=nk$, $v=mk$, $p=n$, $q=m$ yields the following relation
\begin{eqnarray*}
\delta(m+n,k)&\geq& \log c+\delta(n,k)+\log\mathbb P(T_{mk}\leq m-r, T_{mk}\leq D_1)
\\&\geq&  \log c+\delta(n,k)+\log\mathbb P(T_{mk+Lr}\leq m, T_{mk+Lr}\leq D_1)\\
&\geq&\log c+\delta(n,k)+\delta(m,k'),
\end{eqnarray*} 
where $k'$ is any real number satisfying $mk'\geq mk+Lr$. From now on we will write $c$ 
instead of $\log c$. In other words, for each $k$, and each $k'>k$ we have 
\begin{eqnarray}\label{ineq_delta}
\delta(m+n,k)&\geq &c+\delta (n,k)+\delta(m,k'),\end{eqnarray}
for all $m$, $n$ such that $m\geq \frac{Lr}{k'-k}$.
Let $$\underline{\delta}(k)=\liminf_{n\rightarrow\infty} \frac{\delta(n,k)}n,\;\;\;
\overline{\delta}(k)=\limsup_{n\rightarrow\infty}\frac{\delta(n,k)}n.$$

If $k<k'$ then for each $\alpha<\overline \delta(k')$ there exists a sequence 
$m_t$ that goes to infinity as $t\to\infty$ and that satisfies 
$\frac{\delta (m_t,k')}{m_t}\geq \alpha$.
Let us fix one such $m_t$. Each $n\geq m_t$ can be represented as 
$n=am_t+b$ for some $a\geq 0$ and $b\in\{0,1,2,\dots, m_t-1\}$. 
Then we have 
\begin{eqnarray*}\frac{\delta (n,k)}n&=&
\frac{
\delta(am_t+b,k)}{am_t+b}\geq \frac{a\delta(m_t,k') +\delta(b,k)+ac}{am_t+b}\\
&\geq&
\frac{\delta(m_t,k')}{m_t}\cdot\frac1{1+\frac{m_t-1}{am_t}}+
\frac{\delta(b,k)}{am_t+b}+\frac{c}{m_t+\frac ba}\\
&\geq& \frac{\alpha}{1+\frac{m_t-1}{am_t}}+\frac{b,k}{am_t+b}+\frac{c}{m_t+\frac{b}a}
.
\end{eqnarray*}
Keeping $m_t$ fixed and letting $a\to\infty$ we get that $\underline{\delta}(k)\geq \alpha$ 
for all $\alpha<\overline\delta(k')$. This immediately yields  
$\underline \delta(k)\geq \overline \delta(k')$ for all $k<k'$. 
Obviously $\overline\delta(k)\geq \underline\delta(k)$ hence we have 
$$\underline\delta(k')\leq \overline\delta(k')\leq \underline\delta(k)\leq 
\overline\delta(k),\;\;\mbox{for all } k<k'.$$ 
A consequence of the previous inequality is the monotonicity of the functions 
$\overline\delta$ and $\underline \delta$. They are both non-increasing. 
Let $\alpha$ and $\beta$ be two positive rational numbers such that 
$\alpha+\beta=1$. Let $k_1$ and $k_2$ be any two positive real numbers.  According to (\ref{sub_add}) we know that  for each $n\in\mathbb N$ we have:
\begin{eqnarray*}\delta(n,\alpha k_1+\beta k_2)&=&\log\mathbb P(T_{n(\alpha k_1+\beta k_2)}\leq \alpha n+\beta n, T_{n(\alpha k_1+\beta k_2)}\leq D_1)\\&\geq&\log\mathbb P(T_{n\alpha k_1}\leq \alpha n, T_{n\alpha k_1}\leq D_1)+ \\&&
\log\mathbb P(T_{n\beta k_2}\leq \beta n-r, T_{n\beta k_2}\leq D_1)+c\\ &\geq &
\delta(\alpha n, k_1)+\delta (\beta n,k_2')+c,
\end{eqnarray*}
for sufficiently large $n$, where $k_2'$ is any number larger than $k_2$. 
This implies that \begin{eqnarray}\nonumber\overline\delta(\alpha k_1+\beta k_2)&\geq& 
\limsup \frac1n\left(\delta(\alpha n,k_1)+\delta(\beta n,k_2')\right)\\&\geq& 
\alpha\underline\delta(k_1)+\beta\underline\delta(k_2').\label{convexity_rationals}
\end{eqnarray}
Let us justify the second inequality. 
The previous $\limsup$ is definitely larger than the $\liminf$ over the 
sequence of those integers $n$ that are divisible by the denominators of 
both $\alpha$ and $\beta$. 

Consider now two positive real numbers $k<k'$ for which $\underline\delta(k')$ is a real number (i.e. not $-\infty$). Let $\alpha_n$ and $\beta_n$ be two sequences of positive rational numbers such that $\alpha_n\rightarrow 1$, $\beta_n\rightarrow 0$, $\alpha_n+\beta_n=1$. Let $k''$ be a real number such that $k<k''<k'$. Then the inequality (\ref{convexity_rationals}) implies:
\begin{eqnarray}\nonumber\overline\delta(\alpha_nk+\beta_n k'')\geq \alpha_n\underline\delta(k)+\beta_n\underline \delta(k').\end{eqnarray} Taking the limit of both sides as $n\rightarrow\infty$ and using the monotonicity of $\overline\delta$ we get: 
$$\lim_{\varepsilon>0, \varepsilon\rightarrow 0} \overline \delta(k+\varepsilon)\geq \underline \delta(k).$$
This inequality together with  $\underline\delta(k+\varepsilon/2)\geq
\overline\delta(k+\varepsilon)$ implies that $\underline\delta$ is right-continuous. 

Let us now choose the sequences $\alpha_n$, $\beta_n$, $k^1_n$, $k^2_n$, $k^3_n$ that satisfy:
\begin{eqnarray*}
&&\alpha_n,\beta_n\in\mathbb Q_+,\;\; k_n^1,k_n^2,k_n^3\in\mathbb R_+,\\&& \alpha_n+\beta_n=1,\;\; \alpha_n\rightarrow 0,\;\; \beta_n\rightarrow 1,\\&& k^2_n<k^3_n<k<k^1_n, \;\;k^2_n,k^3_n\nearrow k,\;\; k^1_n\searrow k,\\&& \alpha_nk^1_n+\beta_nk^2_n=k.\end{eqnarray*} Then the inequality (\ref{convexity_rationals}) implies  $$\overline\delta(k)\geq \alpha_n\underline\delta(k^1_n)+\beta_n\underline\delta(k^3_n).$$ Letting $n\rightarrow\infty$ and using the monotonicity of $\underline \delta$ gives us:
$$\overline \delta(k)\geq \lim_{\varepsilon>0,\varepsilon\rightarrow0}\underline \delta(k-\varepsilon).$$
Using the fact $\underline \delta(k-\varepsilon)\geq \overline\delta(k-\varepsilon/2)$ gives that $\overline\delta(k)\geq \overline \delta(k-0)$. Here $\overline \delta(k-0)$ is defined as $$\overline{\delta}(k-0)=\lim_{\varepsilon\rightarrow 0+} \overline\delta (k-\varepsilon).$$ In other words, $\overline \delta$ is left-continuous.

Let us now choose the sequences:
\begin{eqnarray*}
&&\alpha_n,\beta_n\in\mathbb Q_+,\;\; k^1_n,k^2_n,k^3_n\in\mathbb R_+,\\&& \alpha_n+\beta_n=1,\;\; \alpha_n\rightarrow 1,\;\; \beta_n\rightarrow 0,\\&& k^3_n>k^2_n>k>k^1_n, \;\;k^2_n,k^3_n\searrow k,\;\; k^1_n\nearrow k,\\&& \alpha_nk^1_n+\beta_nk^2_n>k.\end{eqnarray*}
Placing these sequences in the inequality (\ref{convexity_rationals}) gives us the following relation:
$$\underline\delta(k)\geq\overline\delta(\alpha_nk^1_n+\beta_nk^2_n)\geq \alpha_n\underline\delta(k^1_n)+\beta_n\underline\delta(k^3_n)\geq
\alpha_n\overline\delta(k)+\beta_n\underline\delta(k^3_n).$$
Letting $n\rightarrow\infty$ implies $\underline \delta(k)\geq \overline\delta(k)$. Summing up all the facts we got for $\underline\delta$ and $\overline\delta$ we have that $\underline\delta(k)=\overline\delta(k)=\gamma(k)$ for some function $\gamma$. Moreover, $\gamma$ is continuous and $\gamma(\alpha k_1+\beta k_2)\geq \alpha \gamma(k_1)+\beta\gamma(k_2)$ for $\alpha,\beta\in\mathbb Q_+$ such that $\alpha+\beta =1$. Because of the continuity, the last inequality holds for all $\alpha,\beta\in\mathbb R_+$ such that $\alpha+\beta=1$. This means that $\gamma$ is concave and the equality (\ref{limit_gamma}) is established. 
Now we have the relation (\ref{limit_psi}) and the concavity of the function $\psi$. Namely, Lemma \ref{second_lemma} implies that $\psi=\gamma$. 

Using Lemma \ref{first_lemma} we get that \begin{eqnarray*}\limsup_{n\rightarrow\infty}\frac1n\log\mathbb P(X_n\cdot l\geq nk)&\leq& \psi(k),\;\;\mbox{and}\\
\liminf_{n\rightarrow\infty}\frac1n\log\mathbb P(X_n\cdot l\geq nk)&\geq& 
\psi(k+\varepsilon), \end{eqnarray*}
for all $\varepsilon>0$. If $k$ belongs to the interior of $\psi^{-1}(\mathbb R)$, we can take $\varepsilon\rightarrow0$ in the previous inequality and use the continuity of $\psi$ to obtain $$\liminf_{n\rightarrow\infty}\frac1n\log\mathbb P(X_n\cdot l\geq nk)\geq \psi(k).$$
This in turn implies (\ref{main_limit}) and the concavity of $\phi$. \hfill $\Box$

\vspace{0.5cm}

The Ellis-G\" artner theorem will enable us to get some more information on lower and upper bound large deviations for general sets. 

\begin{definition} Assume that $\Lambda^*$ is a convex conjugate of the function $\Lambda$. We say that $y\in\mathbb R^d$ is an {\em exposed point} of $\Lambda^*$ if for some $\lambda \in\mathbb R^d$ and all $x\neq y$, \begin{eqnarray}\label{dz2.3.4.} \lambda \cdot y- \Lambda^*(y)>\lambda\cdot x- \Lambda^*(x).\end{eqnarray}
$\lambda$ in (\ref{dz2.3.4.}) is called an {\em exposing hyperplane}. (see \cite{dz})
\end{definition}

We are now ready to prove the theorem stated in the introduction.
\begin{theorem}\label{tdz2.3.2.} There exists a convex function $\Lambda:\mathbb R^d\rightarrow\bar{\mathbb R}$ such that $\lim_{n\rightarrow\infty}\frac1n\log\mathbb E\left[e^{X_n\cdot \lambda}
\right]=\Lambda(\lambda)$.
\end{theorem}
\noindent{\bf Proof.} 
As noted in the remark after Theorem \ref{main_theorem} there exists a function $\Phi:\mathbb R^d\rightarrow\bar{\mathbb R}$ such that for all $l\in\mathbb R^d$ and $k\in\mathbb R_+$: $$\lim_{n\rightarrow\infty}\frac1n\log\mathbb P(X_n\cdot l\geq kn)=\Phi\left(\frac1kl\right).$$
For each $\lambda\in\mathbb R^d$ and each $k>0$ we have that \begin{eqnarray*}
\liminf\frac1n\log\mathbb E\left(e^{X_n\cdot \lambda}\right)&\geq&\liminf\frac1n\log\mathbb E\left(e^{X_n\cdot\lambda}\cdot 1(X_n\cdot \lambda>kn)\right)\\&\geq&
\liminf\frac1n\log \mathbb E\left(e^{kn}\cdot 1(X_n\cdot \lambda>kn)\right)=
k+\Phi\left(\frac1k\lambda\right).
\end{eqnarray*}
Moreover, from Theorem \ref{case_k=0} we get \begin{eqnarray*}\liminf\frac1n\log\mathbb E\left(e^{X_n\cdot \lambda}\right)&\geq& \liminf\frac1n\log\mathbb E\left(e^{X_n\cdot \lambda}\right)\\&\geq&
\liminf\frac1n\log\mathbb E\left(e^0\cdot 1(X_n\cdot \lambda\geq 0)\right)=0.
\end{eqnarray*}
Therefore $$\liminf\frac1n\log\mathbb E\left(e^{X_n\cdot \lambda}\right)\geq \max\left\{0, 
\sup_{k>0}\left\{k+\Phi\left(\frac{\lambda}k\right)\right\}\right\}.$$
From the boundedness of jumps of the random walk $X_n$ we have that 
$|X_n\cdot \lambda|<L|\lambda|n$. 
Let $r\in\mathbb N$ and $0=k_0<k_1<k_2<\cdots<k_r=L|\lambda|$. Then we have the equality:
\begin{eqnarray*}
\overline{\lim}\frac1n\log\mathbb E\left(e^{X_n\cdot \lambda}\right)&=&\overline{\lim}\frac1n\log\left(
\mathbb E\left((e^{X_n\cdot \lambda}\cdot 1(X_n\cdot \lambda \leq 0)\right)+
\right.\\&&\left.\sum_{i=0}^{r-1}\mathbb E\left(e^{X_n\cdot \lambda}\cdot 
1(nk_i< X_n\cdot \lambda\leq nk_{i+1})
\right)
\right).\end{eqnarray*}
We now use the standard argument to bound each term of the sum by the maximal term. 
\begin{eqnarray*}\overline{\lim}\frac1n\log\mathbb E\left(e^{X_n\cdot \lambda}\right)
&\leq&\overline{\lim}\frac1n\log\left(\mathbb P(X_n\cdot \lambda\leq 0)+
\sum_{i=0}^{r-1}e^{nk_{i+1}}\cdot \mathbb P(X_n\cdot \lambda>nk_i)\right)\\
&\leq &\overline{\lim}\frac1n\log\left(\mathbb P(X_n\cdot \lambda\leq 0)
%\right.\\&&\left.
+
r\max_{0\leq i\leq r-1}e^{nk_{i+1}}\cdot \mathbb P(X_n\cdot \lambda>nk_i)\right)\\&
=& \max\left\{0,\max_{0\leq i\leq r-1}\left\{k_{i+1}+\Phi\left(\frac{\lambda}{ k_i}\right)\right\}\right\}.
\end{eqnarray*} 
The last equality is true because $r\in\mathbb N$ is a fixed number as 
$n\rightarrow \infty$ and 
Theorem \ref{case_k=0} implies that $\lim_{n\rightarrow\infty}\frac1n\log\mathbb P(X_n\cdot 
(-\lambda)\geq 0)=0$. Theorem \ref{main_theorem} implies that the function $\Phi\left(\frac1k\lambda\right)$ is continuous in $k$ hence taking $r\rightarrow \infty$ and $k_{i+1}-k_i$ constant we get:
$$\limsup\frac1n\log\mathbb E\left(e^{X_n\cdot \lambda}\right)\leq 
\max\left\{0, 
\sup_{k>0}\left\{k+\Phi\left(\frac{\lambda}k\right)\right\}\right\}.$$
This proves the existence of the limit from the statement of the theorem with  $$\Lambda(\lambda)=\max\left\{0, 
\sup_{k>0}\left\{k+\Phi\left(\frac{\lambda}k\right)\right\}\right\}.$$
We will not use this representation for $\Lambda$ to prove its convexity. Notice that all functions $\Lambda_n(\lambda)=\log\mathbb E\left(e^{X_n\cdot \lambda}\right)$ are convex when $n$ is fixed. Indeed, for all $\alpha,\beta\in\mathbb R_+$ with $\alpha+\beta=1$ and all $\lambda,\mu\in\mathbb R^d$ according to the Holder's inequality we have:
\begin{eqnarray*}e^{\Lambda_n(\alpha\lambda+\beta\mu)}&=&\mathbb E\left[\left(e^{X_n\cdot \lambda}\right)^{\alpha}\cdot \left(e^{X_n\cdot\mu}\right)^{\beta}\right]\\&\leq& 
\left(\mathbb E\left[e^{X_n\cdot \lambda}\right]\right)^{\alpha}\cdot 
\left(\mathbb E\left[e^{X_n\cdot \mu}\right]\right)^{\beta}=e^{\alpha \Lambda_n(\lambda)+\beta\Lambda_n(\mu)}.
\end{eqnarray*}
Since the limit of convex functions is convex, as well as the 
maximum of two convex functions, we are able to conclude that $\Lambda$ is convex. 
\hfill $\Box$

\vspace{0.5cm}

Since $\Lambda(\lambda)\leq L|\lambda|$, we have that the origin belongs to the interior 
of the set 
$\{\lambda\in\mathbb R^d:\Lambda(\lambda)<+\infty\}=\mathbb R^d$. 

\begin{theorem}\label{mainmain_theorem} Let $X_n$ be the previously defined 
deterministic walk in a random environment that satisfies the conditions (i)--(iii). 
Let $\Lambda$ be the function from Theorem \ref{tdz2.3.2.} and let $\Lambda^*$ be 
its convex conjugate. Let $\mathcal F$ be the set of exposed points of $\Lambda^*$.
For any closed set $F\subseteq \mathbb R^d$, 
$$\limsup \frac1n \log\mathbb P\left(\frac1nX_n\in F\right)\leq -\inf_{x\in F} 
\Lambda^*(x),$$ 
and for any open set $G\subseteq \mathbb R^d$
$$\liminf \frac1n\log\mathbb P\left(\frac1nX_n\in G\right)\geq 
-\inf_{x\in G\cap \mathcal F} \Lambda^*(x).$$
Moreover, there exists $\delta>0$ such that $\Lambda^*(x)<+\infty$ for $|x|<\delta$. 
\end{theorem}
\noindent{\bf Proof.} 
The conditions for Ellis-G\" artner theorem are now satisfied because of Theorem \ref{tdz2.3.2.}. Direct application of that result proves the first part of the statement. 

For the second part, we will use Lemma \ref{uniformity_nice}. There exists $\kappa>0$ such that for each $\lambda\in\mathbb R^d$ there exists $i\in\{1,2,\dots, m\}$ such that $\lambda\cdot u_i>\kappa|\lambda|$.
Then we have $\Lambda(\lambda)\geq \kappa|\lambda|+\Phi\left(\frac{\lambda}{\kappa|\lambda|}\right)$ and:
\begin{eqnarray*}
\Phi\left(\frac{\lambda}{\kappa|\lambda|}\right)&\geq&\liminf \frac1n\log\mathbb P\left(\xi_1=\xi_2=\cdots=\xi_n=u_i\right)\geq c,\end{eqnarray*}
for some constant $c$. Therefore $\Lambda(\lambda)\geq \kappa|\lambda|+c$, and 
$$\Lambda^*(x) \leq \sup_{\lambda}\{\lambda\cdot x- \kappa|\lambda|-c\} \leq 
-c+\sup_{\lambda}\{|\lambda|\cdot |x|-\kappa|\lambda|\}.$$ 
Hence if $|x|\leq \kappa$ then $\Lambda^*(x)\leq -c<+\infty$.   \hfill $\Box$

\section{Law of Large Numbers} \label{section_lln}
Let us conclude this study of deterministic walks with some notes about the law of large numbers. We will prove that our walk will have $0$ as the limiting velocity. We will also prove that the walk will almost surely end in a loop.

\begin{theorem}\label{lln}
If $X_n$ is defined as before then $$\lim_{n\rightarrow\infty} \frac1n\mathbb E(X_n)=0.$$
\end{theorem}

\noindent{\bf Proof.} It suffices to prove that $\lim_{n\rightarrow\infty} \frac1n\mathbb E(X_n\cdot l)=0$ for each $l\in\mathbb R^d$, because the zero vector is the only one orthogonal to the entire $\mathbb R^d$. Furthermore, the problem can be reduced to proving that $\frac1n\mathbb E[X_n\cdot l]^+$ converges to $0$ because $X_n\cdot l=(X_n\cdot l)^++(X_n\cdot (-l))^+$. By Fubini's theorem we have 
$$\mathbb E\left[\left(\frac1nX_n\cdot l\right)^+\right]=\int_0^{+\infty} \mathbb P(X_n\cdot l>nt)\,dt. $$
Since $\{X_n\cdot l>nt\}=\emptyset$ for $t>L$ the previous integration could be performed on the interval $(0,L)$ only. Let $x_1$, $\dots$, $x_s$ be a sequence from Theorem \ref{theorem_1}, and let $y_k=\sum_{i=1}^k x_i$. Define the random walk $Y_i$ as $Y_i=X_{s+i}$. 
The probability that the walk will reach the half-space  $H_{nt}^l$ before time $n$ is smaller than the probability of the following event: The walk does not make a loop in first $s$ steps, and after that it reaches the half-space $H_{nt-sL}^l$. Therefore we deduce that for each $t\in(0,L)$ the following inequality holds: 
\begin{eqnarray*}
\mathbb P(X_n\cdot l>nt)&\leq&\mathbb E\left[1((X_1,\dots, X_s)\neq (y_1, \dots, y_s))\cdot 
1(Y_{n-s}\cdot l\geq nt-sL)\right]\\
&=&\mathbb E\left[\mathbb E\left[1((X_1,\dots, X_s)\neq (y_1, \dots, y_s))
%\right.\right.\\&&\left.\left.
\cdot 
1(Y_{n-s}\cdot l\geq nt-sL)|\mathcal F_{Y_1, \dots, Y_{n-s}}
\right]\right],
\end{eqnarray*}
where $\mathcal F_{Y_1, \dots, Y_s}$ is a $\sigma$-field determined by $\eta_z$ for all
$z\in\mathbb Z^d\setminus\{X_1, \dots, X_s\}$ such that $\min_{i=1}^n|z-X_i|\leq M$. 
The previous inequality now implies that 
\begin{eqnarray*}
\mathbb P(X_n\cdot l\geq nt)&\leq & \mathbb E \left[1(Y_{n-s}\geq nt-sL)\cdot\right.\\
&&\left. \mathbb E\left[1((X_1,\dots, X_s)\neq (y_1, \dots, y_s))|\mathcal F_{Y_1,\dots, Y_{n-s}}\right]\right].
\end{eqnarray*} 
From Theorem \ref{theorem_1} we have that $\mathbb E\left[1((X_1,\dots, X_s)\neq (y_1, \dots, y_s))|\mathcal F_{Y_1,\dots, Y_{n-s}}\right]\leq 1-c$ for some constant $c>0$. Let us denote $g=1-c$. We know that $g\in(0,1)$. Using mathematical induction, we can repeat the previous sequence of inequalities $[nt/sL]$ times to obtain that $\mathbb P(X_n\cdot l\geq nt)\leq g^{[nt/sL]}$. 
Now we have that for all $t_0>0$ the following inequality holds:
\begin{eqnarray*}\frac1n\mathbb E(X_n\cdot l)^+&=&\int_0^L\mathbb P(X_n\cdot l>nt)\,dt\\&=&\int_0^{t_0}\mathbb P(X_n\cdot l>nt)\,dt+\int_{t_0}^L\mathbb P(X_n\cdot l>nt)\,dt\\&\leq & t_0+(L-t_0)\cdot g^{[nt_0/sL]}\\&\leq& t_0+L\cdot g^{[nt_0/sL]}.
\end{eqnarray*}
If we keep $t_0$ fixed and let $n\rightarrow \infty$ it is easy to see that the last quantity converges to $0$. 
Therefore $\limsup \frac1n \mathbb E(X_n\cdot l)^+\leq t_0$. However, this holds for every $t_0>0$ hence 
$\limsup \frac1n\mathbb E(X_n\cdot l)^+\leq 0$. This finishes the proof of the theorem.
\hfill $\Box$

\vspace{0.5cm}

\begin{theorem}\label{existence_of_loop_prob1}
 Under the stated assumptions  the walk makes a loop with 
probability $1$.
\end{theorem}
\noindent{\bf Proof.}
Let $x_1$, $\dots$, $x_s$ be the sequence obtained in Theorem \ref{theorem_1}. Define the 
sequence $y_1$, $\dots$, $y_s$ as: $y_i=x_1+\cdots+x_{i}$ for $i\in\{1,2,\dots, s\}$. The sequence $y_1$, $y_2$, $\dots$, $y_s=0$, is a loop.  Theorem \ref{theorem_1} required us to fix a vector $l$. We are guaranteed that the sequence $y_1$, $\dots$, $y_s$ will be in the half-space orthogonal to $l$. It really does not matter which $l$ was chosen, since any loop will be sufficient for the argument that follows.  

Denote by $B_n$ 
the event that $X_{ns+j}\in\{X_1,\dots, X_{ns+j-1}\}$ for some $j\in\{1,2,\dots, s\}$. 
Let $N$ denote the event that the walk doesn't have loops. Then for each $n\in\mathbb N$ we have  
$\mathbb P(N)\leq \mathbb P(B_0^C\cap B_1^C\cap\cdots\cap B_n^C)$. 
We will prove that the latter probability converges to $0$ as $n\to \infty$. 
We start by writing
 
\begin{eqnarray*}\mathbb P\left(B_0^C\cap B_1^C\cap\cdots \cap B_n^C\right)&=&
\mathbb P\left(B_0^C\cap B_1^C\cap\cdots \cap B_{n-1}^C\right)-\mathbb P
\left(B_0^C\cap B_1^C\cap\cdots \cap B_{n-1}^C\cap B_n\right).
\end{eqnarray*}
For $i\in\{1,2,\dots, s\}$ consider the following events: 
$$K_i =\left(\bigcap_{j=0}^{i-1}\left\{X_{ns+j}\not\in\{X_1,\dots, X_{ns+j-1}\}\right\}\right)\cap 
\left\{X_{ns+i}\in\{X_1,\dots, X_{ns+i-1}\}\right\}.$$ 
The events $K_i$ are all disjoint and $B_n=\bigcup_{i=1}^s K_i$ hence 
\begin{eqnarray*}
\mathbb P
\left(B_0^C\cap B_1^C\cap\cdots \cap B_{n-1}^C\cap B_n\right)&=&
\mathbb E\left[1(B_0^C)\cdot 1(B_1^C)\cdots 1(B_{n-1}^C)\cdot 1(B_n)\right]\\
&=& \sum_{i=1}^s \mathbb E\left[1(B_0^C)\cdot  1(B_1^C)\cdots 1(B_{n-1}^C)\cdot 
1(K_i)\right].
\end{eqnarray*} 
Consider the events:
\begin{eqnarray*}
 K_i'& =&\left\{\{X_{ns}+y_1,\dots ,X_{ns}+y_{i-1}\}\cap \{X_1,\dots, X_{ns}\}=\emptyset\right\} \cap 
\{X_{ns}+y_i\in\{X_1,\dots, X_{ns} \}\}, \mbox{ and }\\
K_i''&=& \left\{X_{ns+1}=X_{ns}+y_1, X_{ns+2}=X_{ns}+y_2, \dots, X_{ns+i}=X_{ns}+y_i
\right\}.
\end{eqnarray*}
We first notice that $K_i'\cap K_i''\subseteq K_i$. The events $K_i'$ are all disjoint and $K_s'=\left(\bigcup_{i=1}^{s-1} K_i'\right)^C$ which implies that
$1(K_1')+\cdots+ 1(K_s')=1$. 

Denote by $\mathcal G_n$ the $\sigma$-algebra generated by the environment at the sites $\{X_1, \dots, X_{ns-1}\}$. The sets $B_0$, $B_1$, $\dots$, $B_n$ are all measurable with respect to this $\sigma$-algebra. This is obvious for $B_0$, $\dots$, $B_{n-1}$. The measurability of $B_n$ follows from the following observation: in order to determine whether $B_n$ has happened or not, we actually do not need to know what is the environment at the site $X_{ns}$, because that environment can affect only the next jump, i.e. $X_{ns+1}-X_{ns}$. 
Similarly, for each $i$, the event $K_i'$ is measurable with respect to $\mathcal G_n$ because the sequence $y_1$, $\dots$, $y_s$ is deterministic. 
Therefore
 \begin{eqnarray*}
\mathbb E\left[1(B_0^C)\cdot 1(B_1^C)\cdots 1(B_{n-1}^C)\cdot 
1(K_i)\right]&\geq&
\mathbb E\left[1(B_0^C)\cdot 1(B_1^C)\cdots 1(B_{n-1}^C)\cdot 
1(K_i')\cdot 1(K_i'')\right]
\\
&=&\mathbb E\left[
\mathbb E\left[1(B_0^C)\cdots 1(B_{n-1}^C)\cdot 1(K_i')\cdot 
1(K_i'')\mid \mathcal G_n\right]\right] \\
&=&
\mathbb E\left[1(B_0^C)\cdots 1(B_{n-1}^C)\cdot 1(K_i')\cdot
\mathbb E\left[ 
1(K_i'')\mid \mathcal G_n\right]\right].
\end{eqnarray*}

For each $i\in\{1,2,\dots, s\}$ denote by $\mathcal H_i$ the $\sigma$-algebra generated 
by the environment consisting of the $M$-neighborhood of the set 
$\{X_1,\dots,X_{ns}\}\cup\{X_{ns}+y_1,
\dots,X_{ns}+y_{i-1}\}$ excluding the points $X_{ns}$, $X_{ns}+y_1$, $\dots$, $X_{ns}+y_{i-1}$. 
On the event $B_0^C\cap B_1^C\cap\cdots\cap B_{n-1}^C\cap K_i'$ we have \begin{eqnarray}\label{ineq:K_i''}
\mathbb E\left[1(K_i'')|\mathcal G_n\right]=\mathbb E\left[\mathbb E\left[1(K_i'')|\mathcal H_i\right]|\mathcal G_n\right]\geq c^i\geq c^s.\end{eqnarray}
The inequality (\ref{ineq:K_i''}) is proved using induction on $i$. The proof is exactly the same as the proof of Theorem \ref{theorem_11}.

Summing up inequalities (\ref{ineq:K_i''}) for $i=1,2,\dots, s$ gives us 
$$\mathbb P\left( B_0^C\cap\cdots\cap B_{n-1}^C\cap B_n\right)\geq c^s \mathbb E\left[B_0^C\cap\cdots\cap 
B_{n-1}^C\cdot \sum_{i=1}^s 1(K_i')\right]=
c^s\mathbb P\left(B_0^C\cap\cdots\cap B_{n-1}^C\right).$$ 
Subtracting both sides of the previous inequality from $\mathbb P\left(B_0^C\cap 
\cdots\cap B_{n-1}^C
\right)$ yields
$$\mathbb P(B_0^C\cap\cdots\cap B_{n-1}^C\cap B_n^C)\leq (1-c^s)\cdot \mathbb P(B_0^C\cap\cdots\cap 
B_{n-1}^C).$$
Using induction on $n$ we now prove that 
$$\mathbb P(B_0^C\cap\cdots\cap B_{n-1}^C\cap B_n^C)\leq
(1-c^s)^n.$$ Since $1-c^s<1$, the sequence $(1-c^s)^n$ converges to $0$ 
and we conclude that $\mathbb P(N)=0$. \hfill $\Box$

\section{Processes in Random Environments with Teleports}
\label{section_dpre}
\subsection{Introduction}
Our aim is to study the continuous time process $X_t$ on $(\Omega,\mathbb P)$ that solves the following ODE:
\begin{eqnarray}\label{ode}
\frac{dX_t}{dt}&=&b\left(\tau_{X_t}\omega\right).
\end{eqnarray}

In order to prove Theorem \ref{det_process} we will consider the probabilities 
$\mathbb P(X_t\cdot l>kt)$ for fixed $l\in\mathbb R^d$ and fixed $k>0$. 

Let us first outline the main difficulties we have in trying to implement the proof of Theorem
\ref{main_theorem} to continuous setting. The proof of Lemma \ref{second_lemma} used the fact that by time $n$, the walk could jump only $n$ times over the hyperplane through $0$. The continuous process could jump infinitely many times over that hyperplane and those jumps could happen in relatively short time. Instead, we  will look at a strip around $0$ of positive width. By strip we mean the set $\{x\in\mathbb R^d: x\cdot l\in(\alpha,\beta)\}$
for $\alpha,\beta\in\mathbb R$, and a fixed vector  $l\in\mathbb R^d$.
The process can't travel over the strip infinitely many times, because the speed is finite. Thus we will consider the events that the process does not backtrack over the hyperplane $Z_{-w}$ for suitable $w>0$. 
We will prove that the probability of the event that the process doesn't backtrack over the hyperplane $Z_{-w}$ is comparable to the probability that the process doesn't backtrack over $Z_0$. The difficulty here comes from the fact that the process can approach $Z_0$ with  low speed and it can in some way introduce a lot of dependence in the environment. This can happen especially if  the vector field $b$ is continuous. Similar difficulties are known in the study of a class of processes called Lorentz lattice gases (see, for example, \cite{bunimovich}, \cite{grimmett}). We introduce additional assumptions to have more randomness in the definition of $X_t$ and that randomness helps the process to escape from such environments. So, we will have some assumptions that continuous processes can't satisfy. It turns out that   the bounded process won't satisfy some of the requirements, either.

\subsection{Description of the environment} 
We will assume the existence of teleports and discontinuities. They will
help us build tunnels that can take the process in the directions we wish it to travel. Also we will be able to build traps that can hold the process for long time.

Before giving formal definition of the process and the environment, let us explain the introduction of Poisson teleports and the discontinuities in the environment. This subsection will not be completely rigorous but we believe that the explanations will help in understanding the motivation behind the formal definitions that will follow in Subsection \ref{subsection:definitions_environment}.

We will consider a Poisson point process independent of the environment. We will assume that the balls of radius $r$ centered at the points of the process serve as ``teleports,'' and our process evolves according to (\ref{ode}) until it spends a fixed  positive time in some of the teleports. Once that happens, the process will reappear at another location. In some sense, these teleports correspond to the locations where $b$ is infinite.

In order to use the arguments based on subadditivity we need the stationarity of 
the underlying random field. Our assumptions are going to make the grid $\mathbb Z^d$ 
special, so we can't hope for general stationarity. We will assume only that 
$\mathbb P(b\in A) $ and $\mathbb P(b\in \tau_z A)$ have the same distributions for 
$z\in\mathbb Z^d$ (here $A$ is a subset of the set of bounded piecewise continuous 
vector fields). This is going to be sufficient for our purposes, because we will 
assume that after each teleportation the process will appear at the point that 
experiences the same distribution as the initial point. The initial position of the process is assumed to be chosen uniformly inside the unit cube $[0,1]^d$. The random choice of the point is assumed to be independent from the rest of the environment.

Because of the special role of the grid $\mathbb Z^d$ the vector field $b$ can not be continuous. This will require us to modify the notion of solution to (\ref{ode}).  

The first idea is to associate a random vector field $b_Q$ to each unit cube $Q$ of $\mathbb Z^d$. However, then (\ref{ode}) does not have a unique solution because of the lack of continuity. Thus, for each cube $Q$ of $\mathbb Z^d$ we define a vector field $b_Q$ on an open set that contains $Q$. While being in this larger set, the process follows the ODE, and once it leaves it, it is in the interior of the domain of another vector field. Since these larger sets overlap, some care is needed in defining the solution to (\ref{ode}).
 
We now want to define those larger sets that contain the lattice cubes. This can be done in many ways but here we pick the one that will turn out to be convenient for our later study.

\subsection{Definition of the environment}\label{subsection:definitions_environment}
\begin{figure}
\begin{center}
\includegraphics[width=100mm,height=70mm]{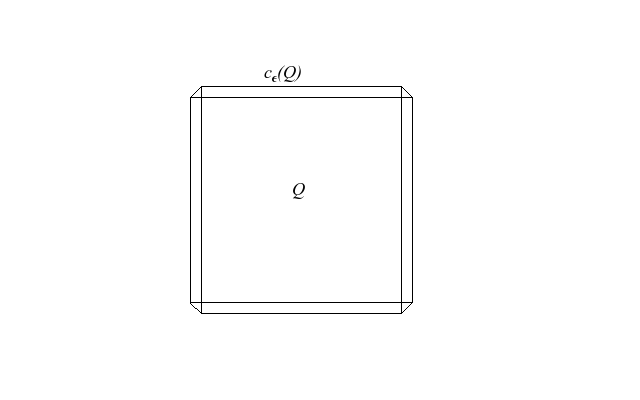} 
\end{center}
\caption{\label{figure_c_varepsilon} The domain $c_{\varepsilon}(Q)$ on which $b_Q$ is defined.
}
\end{figure}

Let $\varepsilon>0$. Consider a unit cube $Q$ whose vertices have all integral coordinates. For each vertex $V$ 
of $Q$ consider the points that are outside of $Q$, belong to the lines containing edges of $Q$, and are 
at a distance $\varepsilon$ from $Q$. There are $d$ such points corresponding to each $V$. Taking 
all $2^d$ vertices of $Q$ we get a total of $d\cdot 2^d$ points.  
Denote by $c_{\varepsilon}(Q)$ the convex hull of all these points. The set 
$c_{\varepsilon}(Q)$ contains 
the cube $Q$ in its interior. Figure \ref{figure_c_varepsilon} shows $c_{\varepsilon}(Q)$ in two dimensions.

Fix $r\in\left(0,\frac1{4\sqrt d}\right)$ and $\varepsilon\in\left(0,\frac r2\right)$. Our environment is given in the following way: 
For each lattice cube $Q$ we have a vector field $b_Q:c_{\varepsilon}(Q)\times \Omega
\to\mathbb R^d$
that is measurable on $c_{\varepsilon}(Q)\times\Omega$.
For each unit cube $Q$  of the lattice, 
denote by $\mathcal F_{Q}$ the $\sigma$-algebra generated by the 
values of $b_{Q'}$ for all cubes $Q'\neq Q$. We make the following assumptions on the vector 
fields $b_Q$:
\begin{enumerate}
\item[(i)] There are positive real numbers $c$ and $\delta_0$ such 
that for each lattice cube $Q$ of edge length  $1$ and each vector 
$l\in\{-1,0,1\}^d\setminus\{0\}$ we have:
$$\mathbb P\left(b_Q(x,\omega)\cdot l>\delta_0 |l| \mbox{ and }\frac{b_Q(x,\omega)\cdot l}{|l|\cdot |b_Q(x,\omega)|}
>\frac{\sqrt{d-1}}{\sqrt d} \mbox{ for all }x\in c_{\varepsilon}(Q)|\mathcal F_{Q}\right)
>c.$$

\item[(ii)]
The vector fields $b_Q$ are bounded, continuous, and  
have finite range dependence, i.e., there are positive constants  $L$ and $M$ such that 
\begin{eqnarray*}
|b_Q(x,\omega)|&\leq& L,\;\;\mbox{ for all }\;\;(x,\omega)\in c_{\varepsilon}(Q)\times\Omega,\\
|b_Q(x,\omega)-b_Q(y,\omega)|&\leq& L|x-y|,\;\;
\mbox{ for all }\;\; x,y\in c_{\varepsilon}(Q)\mbox{ and }\omega\in \Omega,
\end{eqnarray*} 
and $b_Q$ is independent of 
the $\sigma$-field generated by the environment in the cells $Q'$ for which 
 $ \mbox{ dist }(Q',Q)>M$. 
In addition, we make an assumption that the sequence of vector fields $b_Q$ is stationary. 

\item[(iii)] Fix $t_0>0$ such that 
$t_0<\frac{r}{4L}$. Fix also $c_3\in(0,1)$, $\lambda_0>0$, and a sequence 
$(Y_n)_{n=0}^{\infty}$ of iid random variables with values in 
$\{\phi\}\cup [0,1]^d$ that are independent on the environment. 
Each $Y_n$ has a probability $c_3$ of being $\phi$, or (with probability $1-c_3$) 
it is uniformly distributed in $[0,1]^d$.
There exist a set of vectors $\{u_1, \dots, u_m\}\in\mathbb R^d$ 
such that for each $l\in\mathbb R^d\setminus\{0\}$ there is $i\in\{1,2,\dots, m\}$ that satisfies 
$$l\cdot u_i> 3\sqrt d.$$ 
For each $i$, we consider a Poisson point process $\mathcal P_i$ of intensity 
$\lambda_0$. These processes are assumed to be independent among themselves, independent on 
the environment, and on the sequence $(Y_n)$. 
\end{enumerate}

For each point $x$ of the process $\mathcal P_i$, the ball $B_r(x)$ will be called a $(u_i,t_0)$-teleport. In the definitions of the process below, we will see that if the particle spent the time $t_0$ in a $(u_i,t_0)$-teleport, it will have a chance of making a jump in the direction of the vector $u_i$. When there is no ambiguity, a $(u_i,t_0)$-teleport will be called simply $u_i$-teleport. We will often use word teleport for points of Poisson processes (not just for their $r$-neighborhoods). 

\subsection{Description of the process $X_t$}
Before the precise definition of $X_t$ let us describe the process informally. 
The initial position $X_0$ is uniformly distributed inside $[0,1]^d$.

Together with the process we will 
define the sequence of lattice cubes $Q_1=[0,1]^d$, $Q_2$, $\dots$ and the sequence of 
times $\tau_0=0$, $\tau_1$, $\tau_2$, $\dots$ that tell us which $b_Q$ is used 
in the defining ODE.  The time $\tau_n$ can be taught to be  
the ``expiration time'' of the cube $Q_n$, and this is precisely the exit time from the 
interior of the set $c_{\varepsilon}(Q_n)$. The sequence $\tau$ will be increasing, while the cubes in the 
sequence $Q$ may repeat.

The vector field $b_{Q_j}$ is governing the process between times 
$\tau_{j-1}$ and $\tau_j$. Since $b_{Q_j}$ is Lipschitz in $c_{\varepsilon}(Q_j)$ the process $X_t$ is properly defined by the equation (\ref{ode}). Even though these large sets $c_{\varepsilon}(Q)$ intersect, the process does not switch to the new vector field until it exits the old $c_{\varepsilon}(Q)$. Once it exits $c_{\varepsilon}(Q_j)$, it is going to be in the interior of another cube and then it starts following the ODE corresponding to the vector field of that cube and the initial point that is now in the interior of the domain where $b$ is defined. We may think as of having many vector fields that overlap, but the process follows one vector field at any given time.

In addition to switching the environments the process may enter a teleport. Teleports are balls of radius $r$ around Poisson points of processes $\mathcal P_i$. Once the process spends a lot of time inside a teleport $i$, it makes a jump. The vector for the jump is $u_i$.

\subsection{Definition of the process $X_t$}
We will define simultaneously the 
sequence $(Q_n)_{n=0}^{\infty}$ of cubes, the sequence $(\tau_n)_{n=0}^{\infty}$ of stopping times, and the 
process $X_t$. 

We start by defining $\tau_0=0$, and $Q_1=[0,1]^d$. The position $X_0$ of the process is chosen uniformly at random inside the cube $[0,1]^d$. 
Assume that $\tau_1$, $\dots$, $\tau_{k-1}$, and $Q_1$, $\dots$, $Q_k$ are defined and that 
the process is defined for $t\in[0,\tau_{k-1}]$. We will now define $Q_{k+1}$, $\tau_k$, 
and $X_t$ for $t\in(\tau_{k-1},\tau_k]$.

The process $X_t$ satisfies \begin{eqnarray}\label{ode_with_Q}\frac{d X_t}{dt}=b_{Q_k}(X_t,\omega)\end{eqnarray} until one of the following 
two 
events happen:
\begin{enumerate}
\item[$1^{\circ}$] There is an integer $i\in\{1,2,\dots, m\}$ for which the process spent the entire time $(t-t_0,t)$ in the union of $(u_i,t_0)$-teleports (notice that teleports may overlap). Then  
at time $t$ the process $X_t$ does one of the following:
\begin{enumerate}
\item[$1.1^{\circ}$] It stays at the same place if $Y_{\lfloor t/t_0\rfloor}=\phi$ or there are two different $i$ and $j$ 
such that the process spent time $(t-t_0, t)$ in both $(u_i,t_0)$ and $(u_j,t_0)$ teleports.
\item[$1.2^{\circ}$] If $Y_{\lfloor t/t_0\rfloor}\neq \phi$, the process $X_t$ will make a jump and reappear in the lattice cube containing   $X_{t}+u_i$ at the relative position 
$Y_{\lfloor t/t_0\rfloor}$ within the cube. 
\end{enumerate} 
In each of the cases $1.1^{\circ}$ or $1.2^{\circ}$ we declare the 
moment $t$ as a moment in which the process is not in any of the teleports.  
(This way the time restarts if the process spent $t_0$ time in a teleport, and we 
won't have infinitely many successive teleportations in a short time.)

If the jump occurred then we define $\tau_k=t$, and $Q_{k+1}$ is the cube 
where the process landed after teleportation. If the jump did not occur the process keeps following (\ref{ode_with_Q}).

\item[$2^{\circ}$] 
The process has hit the boundary of $c_{\varepsilon}(Q_k)$. Then we define $\tau_k=t$ and 
$Q_{k+1}$ is the cube that contains $X_t$. There could be two or more candidates for $Q_{k+1}$ but
this is a zero-probability event, and in this case we could choose in any way. 
\end{enumerate}

Notice that by fixing the time $t_0$ we 
make sure that the average speed of $X_t$ remains bounded. 
Although there could be a teleportation involved, 
the particle has to wait to be teleported and the distance it can go is bounded. 

There could be regions belonging to more than one teleport. If the particle is a subject to two or more different teleportations, then the jump will be suppressed. A particle could be subject to more than one teleportation if it enters a teleport by a jump. Entering two teleports simultaneously by the means of diffusion is a zero-probability event. 

The requirement $t_0<\frac{r}{4L}$ guarantees that if the process comes to within $r/2$ of a Poisson point, it is going to be teleported for sure, because it won't have enough time to escape. 

In the rest of the paper we will write just $b$ instead of $b_Q$ whenever it is clear 
from the context which $Q$ we are dealing with. 

\subsection{Main Results}
We consider the moment-generating 
function $\mathbb E(\exp(\lambda\cdot X_t))$ for $\lambda\in\mathbb R^d$. Our goal is to prove that there is a convex function $\Lambda$ such that
$$\lim_{t\to+\infty}\frac1t\log\mathbb E(\exp(\lambda\cdot X_t))=\Lambda(\lambda).$$ 

We now fix a vector $l\in\mathbb R^d\setminus\{0\}$ and recall the definition of  the hitting times $T_p$ (for $p\in\mathbb R_+$) of hyperplanes  from Section \ref{hitting_times}. 
%Let us point out that condition (i) gives us that there is a positive constant $\delta_0$ (we call it $\delta_0$ again as we may decrease the original $\delta_0$ if necessary) such that $\mathbb P(b(x,\omega)\cdot l>\delta_0|l|)>c$.

We have an analogous result to Lemma \ref{first_lemma}.

 \begin{lemma}\label{first_lemma_cont}
The following inequality holds: 
\begin{eqnarray*}
\limsup_{t\rightarrow +\infty}\frac1t \log \mathbb P(X_t\cdot l\geq tk)&\leq&
\limsup_{t\rightarrow +\infty} \frac1t \log\mathbb P(T_{tk}\leq t)
\end{eqnarray*}
In addition, for each $\varepsilon>0$ we have:
\begin{eqnarray*}
\liminf_{t\rightarrow+\infty} \frac1t \log\mathbb P(T_{t(k+\varepsilon)}\leq t)&\leq&
\liminf_{t\rightarrow +\infty}\frac1t \log \mathbb P(X_t\cdot l\geq tk).
\end{eqnarray*}
\end{lemma} 

\noindent{\bf Proof.}
The first inequality follows immediately as in the discrete case. For the second one we have to modify the argument a bit. We are not able to construct a loop as we did in the discrete case. The reason is that the curve $X_t$ may not have a self intersection. It is known that in the case of continuous gradient vector fields the process will never have a loop.

Let us define $$S=\{(x_1,\dots, x_d)\in\mathbb R^d: x_i\in\{-1,0,1\}, i=1,2,\dots, d\}. $$
Denote the points of $S$ by 
$P_0$, $P_1$, $P_2$, $\dots$, $P_{3^d-1}$. Assume that $P_0$ coincides with the 
origin $O$. Let $C_i$ ($0\leq i\leq 3^d-1$) be the translation of the cube $[0,1]^d$ for the vector $\overrightarrow{OP_i}$. 
For each integer $i\in\{1,2,\dots, 3^d-1\}$, consider the event 
\begin{eqnarray}D_i=\left\{\omega\in\Omega: 
b_{C_i}(z,\omega)\cdot \overrightarrow{P_iO}>\frac{\delta_0}2|\overrightarrow{P_iO}| \mbox{ and }
\frac{b_{C_i}(z,\omega)\cdot \overrightarrow{P_iO}}
{ |b_{C_i}(z,\omega)|\cdot |\overrightarrow{P_iO}|}>
\frac{\sqrt{d-1}}{\sqrt d} \mbox{ for all } z\in c_{\varepsilon}(C_i)
\right\}.\label{environment_D_i}\end{eqnarray}
The constant $\delta_0>0$ is chosen according to the assumption (i) on the vector fields $b_Q$. Let $D=\bigcap_{i=1}^{3^d-1}D_i\cap D'$, where $D'$ is the event that there are no 
Poisson points in the $r$-neighborhood of $c_{\varepsilon}(C_0)$. Denote by $\mathcal F_{C_1,\dots, C_{3^d-1}}$ the sigma algebra generated by the environment influencing $b_{Q'}$ for $Q'\not\in\{C_1,\dots, C_{3^d-1}\}$. According to our assumptions, we know that there is a constant $c'>0$ such that $$\mathbb P\left(D\mid\mathcal F_{C_1,\dots, C_{3^d-1}}\right)>c'.$$ 

We now prove that if the process ever enters the cube $C_0$ and gets under the influence of $b_{C_0}$, it will stay forever in the closure of $c_{\varepsilon}(C_0)$. We will show that if $Q_k=C_0$, then for each $t\geq t_{k-1}$ we have $X_t\in \overline{c_{\varepsilon}(C_0)}$. Assume the contrary, that $X_t\in C_i\setminus \overline{c_{\varepsilon}(C_0)}$ for some $i\in\{1$, $\dots$, $3^d-1\}$. We may also assume that $C_i$ is the first cube to which $X_t$ has escaped. Let $Y_t=\overrightarrow{OX_t}\cdot \overrightarrow{P_iO}$. Set $$t_1=\inf\{s: X_u \in C_i\setminus c_{\varepsilon}(C_0) \mbox{ for all } u\in(s,t)\}.$$ Then $X_{t_1}$ belongs to the boundary of $c_{\varepsilon}(C_0)$ and because of our choice 
of $c_{\varepsilon}$ we can guarantee that $Y_t<Y_{t_1}$. Using the fundamental theorem of calculus we obtain $$0>Y_t-Y_{t_1}=\int_{t_1}^t  \overrightarrow{P_iO}\cdot b(X_s,\omega)\,ds>0,$$ since $\overrightarrow{P_iO}\cdot b(X_s,\omega)>0$ for all $s\in(t_1,t)$ on the event $D_i$. This is a contradiction.

Now we can finish the proof in the similar way as in the discrete case. 
Let   $S_{kt}$ denote the ``shard''-like surface consisting of faces of the
grid of size $3$ that is in front of the plane $Z_{tk}$ (here ``in front of'' 
means with respect to the direction $l$). 
Here is the precise definition of $S_{tk}$. For each point $x$ of the 
hyperplane $Z_{tk}$, consider all points $z$ of the integral lattice 
such that $z\cdot l>x\cdot l$ and denote by $z_x$ the one that minimizes 
the distance from the point $x$. Consider the hyperplanes of the integral 
lattice that pass through $z_x$. They divide the space in $2^d$ subsets 
($d$-dimensional dihedral angles) only one of which does not contain points of $Z_{tk}$ in its interior. 
For every $x$ we form one such subset (notice that there will be only countably 
many of them). 
Denote by $\Sigma_{tk}$ their union, and $S_{tk}$ is the boundary 
\begin{figure}
\begin{center}
\includegraphics[width=70mm,height=70mm]{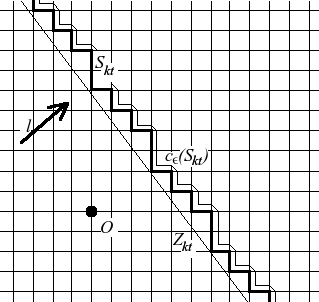} 
\end{center}
\caption{\label{figure_shard} Surfaces $c_{\varepsilon}(S_{kt})$, $S_{kt}$, hyperplane $Z_{kt}$,  and  vector $l$}
\end{figure}
of $\Sigma_{tk}$. The Figure 
\ref{figure_shard} shows how $S_{tk}$ looks in $2$ dimensions.  
Let $$\hat T_{tk}=T_{c_{\varepsilon}(S_{tk})},$$ where $c_{\varepsilon}(S_{tk})$ is defined in the 
following way: Denote by $\mathcal Q$ the set of all lattice cubes $Q$ that share points with $S_{tk}$ but are outside of 
$\Sigma_{tk}$. Denote  $$U_{tk}=\bigcup_{Q\in\mathcal Q}c_{\varepsilon}(Q),$$ and take 
$c_{\varepsilon}(S_{tk})$ to be the portion of the boundary of $U_{tk}$ that is inside 
$\Sigma_{tk}$. 

Denote by $C$ the translation of cube $\bigcup_{i=0}^{3^d-1} C_i$ with the following properties: $C$ contains the point $X_{\hat T_{tk}}$ on one of its faces, and is on different side of $S_{tk}$ than the origin. Denote by $\hat D$ the event that the environment in $C$ is the translation for of the environment defined by (\ref{environment_D_i}) and the paragraph below it. Denote by $\mathcal F_{\hat T_{kt}}$ the $\sigma$-algebra 
of the environment and teleports determined by the cubes that the walk was influenced by the time $\hat T_{kt}$. 
Using conditional expectations we get
 \begin{eqnarray} \nonumber
\mathbb P(X_t\cdot l\geq kt)&=& \mathbb P(X_t\cdot l\geq kt, T_{kt}\leq t)\\
\nonumber &\geq&
\mathbb P\left( X_t\cdot l \geq kt, \hat T_{kt}\leq t, \hat D\right)
\\ \nonumber
&
=
&\mathbb E\left[\mathbb E\left[
1(\hat T_{kt}\leq t)\cdot 1(X_t\cdot l\geq kt)
\cdot 1(\hat D)\mid \mathcal F_{\hat T_{kt}}\right]\right]
\\ \nonumber &=&
\mathbb E\left[
1(\hat T_{kt}\leq t)\cdot
\mathbb E\left[ 1(X_t\cdot l\geq kt)
\cdot 1(\hat D)\mid \mathcal F_{\hat T_{kt}}\right]\right].
\end{eqnarray}
Denote by $R$ the right-hand side of the last inequality. 
Let $\hat F$ be the event that there exists $i\in\{1,\dots, m\}$ and 
$t_1\in \left(t_0\lfloor \frac{\hat T_{tk}}{t_0}\rfloor, 
\hat T_{tk}
\right)$ such that $X_t$ spent all the time $(t_1-t_0,t_1)$ in a $u_i$-teleport. 
Then we can write the quantity $R$ in the following way:  
\begin{eqnarray}\label{use_of_time}
R&=&\mathbb E\left[
1(\hat T_{kt}\leq t)\cdot 1(\hat F)\cdot
\mathbb E\left[ 1(X_t\cdot l\geq kt)
\cdot 1(\hat D)\mid \mathcal F_{\hat T_{kt}}\right]\right]\\
\nonumber
&&+
\mathbb E\left[
1(\hat T_{kt}\leq t)\cdot 1(\hat F^C)\cdot
\mathbb E\left[1(X_t\cdot l\geq kt)
\cdot  1(\hat D)\mid \mathcal F_{\hat T_{kt}}\right]\right].
\end{eqnarray}
Consider the first summand on the right hand side of the previous inequality. 
On the event $\hat F$ we know that 
the random variable $Y_{\lfloor \hat T_{tk}/t_0\rfloor}$ will not be responsible for any further jumps. The process already spent a lot of time in a teleport, so the value 
of $Y_{\lfloor \hat T_{tk}/t_0\rfloor}$ was already used in making the decision whether there will be a jump or not. We don't care about the outcome, because the process reached the level $S_{tk}$, but we know for sure that a single value of $Y$ can't be responsible for two decisions about jumps.
Hence the event 
$$\hat G=\left\{Y_{1+\lfloor \hat T_{tk}/t_0\rfloor}=\cdots=Y_{\lfloor \hat T_{tk}/t_0\rfloor+
\lceil\frac{2\sqrt d}{\delta_0}\rceil}=\phi\right\}$$ is sufficient to assure that the path $X_t$ will enter the trap before having any possible jumps. There are only finitely many terms listed in this sequence, so the probability of  $\hat G$ is strictly positive. 
We also have  that $1(X_t\cdot l\geq kl)=1$ on the intersection 
$\hat D\cap \hat G\cap \hat F\cap \{\hat T_{tk}\leq t\}$. 
Hence there is a constant $c'$ such that on the event $\{\hat T_{kt}\leq t\}\cap\hat F$ we have
$$\mathbb E\left[1(X_t\cdot l\geq kt)\cdot 1(\hat D)\mid \mathcal F_{\hat T_{kt}}\right]
\geq \mathbb E\left[1(X_t\cdot l\geq kt)\cdot 1(\hat D)\cdot 1(\hat G)\mid \mathcal F_{\hat T_{kt}}
\right]=\mathbb E\left[ 1(\hat D)\cdot 1(\hat G)\mid \mathcal F_{\hat T_{kt}}
\right]
\geq c'. $$
Thus we can bound  from below the first summand on the right of (\ref{use_of_time}) by 
$c'\cdot \mathbb P(\hat T_{kt}\leq t, 
\hat F)$.

Let us now consider the second summand. Denote by $\hat E$ the event that
$Y_{\lfloor\hat T_{tk}/t_0\rfloor}=\phi$.  Notice that
\begin{align}\nonumber
\mathbb E\Big[
1(\hat T_{kt}\leq t)& 1(\hat F^C)\cdot
\mathbb E\left[ \cdot 1(X_t\cdot l\geq kt)
\cdot1(\hat D)\mid \mathcal F_{\hat T_{kt}}\right]\Big]
\\ \label{modify_past}
&\geq
\mathbb E\Big[
1(\hat T_{kt}\leq t)\cdot 1(\hat E)
\cdot 1(\hat F^C)\cdot
\mathbb E\left[1(X_t\cdot l\geq kt)
\cdot  1(\hat D)\cdot 1(\hat G)\mid \mathcal F_{\hat T_{kt}}\right]\Big].
\end{align}
Here $\hat G$ is the same event as above. The identity 
$1(X_t\cdot l\geq kt)=1$ holds as before, because on the set $\hat E\cap \hat G$ the process will enter the trap.  
In the same way as in the case of the first summand of (\ref{use_of_time}) we can bound from below the second expression from
(\ref{modify_past})  by
$c'\mathbb P(\hat T_{kt}\leq t, \hat E, \hat F^C)$.   The unpleasant thing is that by introducing the event $\hat E$ we are imposing conditions on the past. However,
on the event $\{T_{kt}\leq t, \hat F^C\}$ we know that the process reached the level $c_{\varepsilon}(S_{kt})$ by time $t$ and 
it didn't spend sufficient time in a teleport to be considered for a jump 
in which $Y_{\lfloor \hat T_{kt}/t_0\rfloor}$ would play a deciding role. 
Therefore $\hat E$ is independent of $\{\hat T_{kt}\leq t, \hat F^C\}$ hence $$\mathbb P(\hat T_{kt}\leq t, \hat E, \hat F^C)=
\mathbb P(\hat E)\cdot 
\mathbb P(\hat T_{kt}\leq t, \hat F^C)=c_3\cdot
\mathbb P(\hat T_{kt}\leq t, \hat F^C).$$
This allows us to conclude that there is a constant $c''>0$ such that
\begin{eqnarray*}\mathbb P(X_t\cdot l\geq kt)&\geq& c'' \left(\mathbb P(\hat T_{kt}\leq t,\hat F)+
\mathbb P(\hat T_{kt}\leq t,\hat F^C)\right)=c'' \mathbb P(\hat T_{kt}\leq t)\\
&\geq& c''\mathbb P(T_{kt+5\sqrt d}\leq t)
.\end{eqnarray*}
Taking the logarithm of both sides of the last inequality, dividing by $t$, and 
taking the $\liminf$ as $t\to+\infty$ one obtains 
the following inequality $$\liminf_{t\to+\infty}
\frac1t\log\mathbb P(X_t\cdot l\geq kt)\geq \liminf_{t\to+\infty}
\frac1t\log\mathbb P(T_{kt+5\sqrt d}\leq t). $$
For each $\varepsilon$, there exists $\tilde t$ such that every $t>\tilde t$ satisfies 
$\left\{T_{(k+\varepsilon)(t)}\leq t-s\right\}\subseteq\left\{T_{kt+5\sqrt d}\leq t-s\right\}$.
Thus $$\liminf_{t\to+\infty}
\frac1t\log\mathbb P(X_t\cdot l\geq kt)\geq
\liminf_{t\to+\infty}
\frac1t\log\mathbb P(T_{(k+\varepsilon)t}\leq t).$$
This completes the proof of the lemma. \hfill $\Box$

\vspace{0.5cm}

Following the approach from the discrete case our goal is to establish a statement similar to Lemma \ref{second_lemma}. Let us define $D_1=\inf\{t>0: X_t\cdot l\leq 0\}$. This is analogous to a stopping time we used in the discrete setting: time of first return to $0$.

\begin{lemma}\label{second_lemma_cont}
Let $k$ and $k'$ be two real numbers such that $0<k'<k$. Then the following two inequalities hold: 
\begin{eqnarray*}\limsup \frac1t\log\mathbb P(T_{tk}\leq t) &\leq& \limsup\frac1t
\log\mathbb P(T_{tk'}\leq t, T_{tk'}\leq D_1)\\
\liminf \frac1t\log\mathbb P\left(T_{tk}\leq t, T_{tk}\leq D_1\right)&\leq &\liminf
\frac1t\log\mathbb P\left(T_{tk}\leq t\right).\end{eqnarray*}
\end{lemma}  
\noindent{\bf Proof.} 
The second inequality is obvious since
$\{T_{kt}\leq t, T_{kt}\leq D_1\}\subseteq \{T_{kt}\leq t\}$. 
Our proof of the first inequality from Lemma \ref{second_lemma} used the fact that in finite time the walk can make only finitely many crossings over the fixed hyperplane. Obviously, we can't use that fact in the continuous setting. The idea is 
to break the process between its crossings of a strip between hyperplanes
$Z_{-w}$ and $Z_0$, where 
$w$ is some fixed real number from the interval $(\frac12,1)$.
Define the following stopping times:
$G_0=0$, $F_0=\inf\{t: X_t\cdot l\leq -w\}$. Having defined $G_i$ and $F_i$, for $i\geq 0$, we inductively define 
\begin{eqnarray*}
G_{i+1}&=&\inf\{t\geq F_i: X_t\cdot l\geq 0\}\\
F_{i+1}&=&\inf\{t\geq G_{i+1}:X_t\cdot l\leq -w\}.
\end{eqnarray*}

We will need the following lemma:
\begin{lemma}\label{t_lemma}
For any two real numbers $k$ and $k'$ satisfying 
$0<k'<k$ we have 
$$\limsup \frac1t \log\mathbb P(T_{tk\leq t}, T_{tk}\leq F_0)\leq 
\limsup \frac1t\log\mathbb P(T_{tk'}\leq t, T_{tk'}\leq D_1).$$
\end{lemma}
\noindent{\bf Proof.}
There is at least one of the vectors from the assumption (iii), say $u_1$, such that 
$l\cdot u_1>3\sqrt d$. 
Let $\Gamma$ be the event that every point in the cube $[0,1]^d$ is at a distance smaller than $\frac12r$ of a $(u_1,t_0)$ teleport, every point in the cube $[0,1]$ is at a distance at least $r$ from any other teleport, and that $Y_0\neq\phi$. The conditional probabilities of that event are bounded below by a constant $c$. By the time $t_0$ the process will be away from the origin. Denote by $Q$ the cube to which it gets teleported. The cube $Q$ is the translation of $[0,1]^d$ for some vector $q\in\mathbb Z^d$. Denote by $\tilde O$, $\tilde Z_0$, and $\tilde Z_{-u_1\cdot l}$ the translations of the origin $O$, hyperplane $Z_0$, and hyperplane $Z_{-u_1\cdot l}$ by the vector $q$. 
Denote by $\tilde X$ the process defined by 
$\tilde X_t=X_{t+t_0}$.  Let $\tilde G_i$ and $\tilde F_i$ denote the stopping times corresponding to the process $\tilde X$ which are analogous to the stopping times 
$G_i$ and $F_i$. 
Consider the projection $\tilde O^l$ of the point $\tilde O$ to the hyperplane $\tilde Z_{-u_1\cdot l}$. Consider the $d-1$-dimensional ball $B_{\max_i\{|u_i|\}}(O^l)$ in $\tilde Z_{-u_i\cdot l}$, and let $B'$ be the $1$-neighborhood of this ball in the space $\mathbb R^d$.
Let $\tilde \Gamma$ denote the event that each point in $B'$ is at a distance smaller than $\frac 12r$ to a $(u_1,t_0)$ teleport and larger than $r$ from any other teleport.
The probability of that event is strictly positive and the event is independent on 
the $\sigma$-algebra generated by the environment in the positive  $l$-direction of the hyperplane $\tilde Z_0$.
Therefore we have 
\begin{eqnarray*}
\mathbb P(T_{tk'}\leq t, T_{tk'}\leq D_1)&\geq&
\mathbb P(T_{tk'}\leq t, T_{tk'}\leq D_1, \Gamma)\\
&\geq &\mathbb P(\tilde T_{tk'}\leq t-t_0, \tilde T_{tk'}\leq \tilde F_0, \tilde \Gamma)\\
&\geq& c\cdot\mathbb P(T_{tk'}\leq t-t_0, T_{tk'}\leq F_0).
\end{eqnarray*}
It remains to notice that for sufficiently large $t$ we have that 
$tk'\leq (t-t_0)k$, hence $$\mathbb P(T_{tk'}\leq t-t_0)\geq 
\mathbb P(T_{(t-t_0)k}\leq t-t_0). $$
This completes the proof of Lemma
\ref{t_lemma}. \hfill $\Box$

\vspace{0.5cm}

We now return to the proof of Lemma \ref{second_lemma_cont}. 
 For any real number $u>0$  we write the event $\{T_{u}\leq t\}$ as the following union:
$$\{T_{u}\leq t\}=\{T_{u}\leq t, T_{u}\leq F_0\}\cup \bigcup_{i=0}^{\infty}
\{T_{u}\leq t, F_i\leq T_{u}\leq F_{i+1}\}. $$ The last union turns out to be finite, because we can prove that if $T_u\leq t$, then $T_u\leq F_{\lceil tL\rceil}$. 
We first find a lower bound on $F_{\lceil u\rceil}$. 
The fundamental theorem of calculus implies
$$-w=X_{F_i}\cdot l-X_{G_i}\cdot l=\int_{G_i}^{F_i}b(\tau_{X_s}\omega)\cdot l\,ds
\geq -L\cdot |l| \cdot (F_i-G_i)$$ which together with $|l|=1$ yields 
$F_i-G_i\geq \frac wL$. Using this inequality we obtain
\begin{eqnarray*}
F_{m}&=&F_0+\sum_{i=1}^{m} (F_i-F_{i-1}) \geq
\sum_{i=1}^{m}(F_i-G_i) \geq \frac{wm}{L}.
\end{eqnarray*}
Therefore $F_{\lceil tL/w\rceil}\geq t$ and on $\{T_u\leq t\}$ we immediately get 
$F_{\lceil tL/w\rceil}\geq T_u$.  This implies that 
$$\mathbb P(T_u\leq t)\leq \sum_{i=0}^{\left\lceil \frac{tL}{w}\right\rceil}\mathbb P(T_u\leq t, 
F_i\leq T_u\leq F_{i+1}).$$
Let us prove that each term on the right hand side of the last inequality can be bounded by the quantity $\mathbb P(T_u\leq t, T_u\leq F_0)$. 
Denote by $\mathcal F_{G_{i+1},T_u}$ the sigma algebra generated by the environment  
contained in the $M$-neighborhood of the process from $G_{i+1}$ to $T_u$.
Notice that if $T_u\leq t$ and $T_u\geq F_i$, then the process 
has made at least $i$ trips over the region between the hyperplanes 
$Z_0$ and $Z_{-w}$. Since $T_u\leq F_{i+1}$ and $T_u\leq t$ we conclude that the process has crossed the hyperplane $Z_u$ by time $t$ which means that it 
had to cross the hyperplane $Z_0$ again. Therefore $G_{i+1}\leq t$.
Let $\tilde X$ be the process starting at time $G_{i+1}$. More precisely, we define
$\tilde X_t=X_{G_{i+1}+t}$. We use $\tilde F_i$ and $\tilde G_i$ to denote the stopping 
times  for $\tilde X_t$ analogous to $F_i$ and $G_i$.  \begin{eqnarray*}
\mathbb P(T_u\leq t, F_i\leq T_u\leq F_{i+1})&=&\mathbb E\left[
 \mathbb E\left[
1(T_u\leq t) \cdot 1(T_u\geq F_i)\cdot\left. 1(T_u\leq F_{i+1}) 
\right| \mathcal F_{G_{i+1},T_u}
\right]
\right]\\
&=&\mathbb E\left[\mathbb E\left[
1(\tilde T_u\leq t-G_{i+1})\cdot 1(T_u\geq F_i)\cdot 1(\tilde T_u\leq \tilde F_0)
| \mathcal F_{G_{i+1},T_u}
\right]
\right]\\
&\leq&\mathbb E\left[1(\tilde T_u\leq t)\cdot 1(\tilde T_u\leq \tilde F_0)\cdot 
\mathbb E\left[1(T_u\geq F_i)|\mathcal F_{G_{i+1},T_u}
\right]
\right].\end{eqnarray*}
The last conditional expectation can be bounded above by $1$. Although our process does not posses Markov property, we can use the trivial bounds on the 
indicator functions, namely $1(T_u\geq F_i)\leq 1$.
Therefore $$\mathbb P(T_u\leq t)\leq \left\lceil\frac{tL}w\right\rceil \cdot\mathbb P(T_u\leq t, 
T_u\leq F_0).$$ Placing $u=tk$ and using Lemma \ref{t_lemma} we obtain
\begin{eqnarray*}
\limsup\frac1t\log\mathbb P(T_{tk}\leq t)&\leq&\lim_{t\to+\infty} \frac1t\log t+
\limsup\frac1t\log\mathbb P(T_{tk}\leq t, T_{tk}\leq F_0)\\&\leq& \limsup\frac1t\log\mathbb 
P(T_{tk'}\leq t, T_{tk'}\leq D_1),
\end{eqnarray*}
for any two real numbers $k'$ and $k$ such that $0<k'<k$. The proof of Lemma \ref{second_lemma_cont} is complete.\hfill $\Box$

\vspace{0.5cm}

\noindent{\bf Proof of Theorem \ref{det_process}.}
The proof will proceed in the same way as in the discrete case. We will first prove the existence of the concave function $\phi^l:\mathbb R_+\to\bar{\mathbb R}$ such that for all $k\in\mathbb R_+$:
\begin{eqnarray}
\label{cont_probabilities}
\lim_{t\to+\infty}\frac1t\log\mathbb P(X_t\cdot l\geq tk)&=&\phi^l(k).
\end{eqnarray} 
As in the discrete case we see that there is a function $\Phi:\mathbb R^d\to\mathbb R$ such that $\phi^l(k)=\Phi\left(\frac lk\right)$. In order to prove (\ref{cont_probabilities}) we 
will establish the inequality that has precisely the same form as (\ref{sub_add}), with a difference that our stopping times are defined in the continuous setting. The proof will be adapted to the continuous in the similar way as the proof of Lemma \ref{first_lemma_cont}. In the proof of the mentioned lemma we needed to construct a trap that will hold the process. Here we need a tunnel that will take the process to the environment sufficiently far away from the previously visited sites.

We start by considering the parallelepiped $Q=\left[0,\lceil nL\rceil\right]\times [-3,3]^{d-1}$. Here $n$ is an integer that is going to be specified in a moment. It will depend on $d$ and $l$ only.
Similarly to the proof of Lemma \ref{first_lemma_cont}, denote by $S_u$ 
the surface consisting only of faces of the grid of size $\lceil nL\rceil$ that is in front of $Z_u$ when looking from the origin in the direction $l$. Denote by $Q'_u$ the appropriate rotation and translation of $Q$ so 
it ends up on the side of $S_u$ opposite to the origin. Moreover, $Q'_u$ can 
be placed in such a way that the hitting point of $c_{\varepsilon}(S_u)$ by the process $X_t$ is inside $Q_u'$ and located within $\varepsilon$ of the base of $Q'_u$ (here by base we mean a $d-1$-dimensional face congruent to $[-3,3]^{d-1}$).  
We choose the integer $n$ so that there is a hyperplane $Z_y$ orthogonal 
to $l$ that cuts $Q'_u$ in two prisms each of which has the shortest altitude of at 
least $2L$. 
Denote by $Q_u$ the one of the prisms 
that contains the hitting point of the process $X_t$.

\begin{figure}
\begin{center}
\includegraphics[width=100mm,height=70mm]{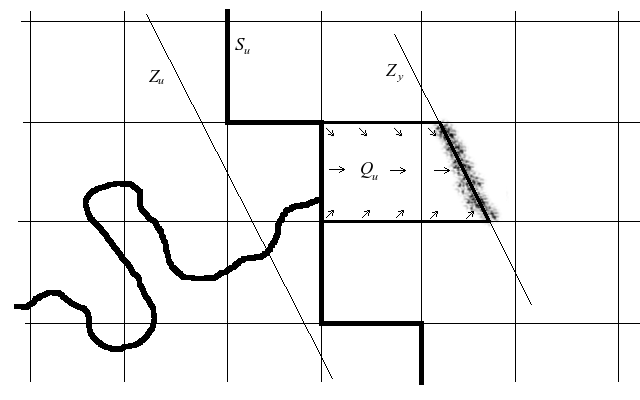} 
\end{center}
\caption{\label{figure_tunnel} The polytope $Q_u$ and teleports near the boundary that is a part of $Z_y$
}
\end{figure}

Denote by $\hat T_u$ the hitting time of $c_{\varepsilon}(S_u)$. 
Denote by $B_1$ the event that after reaching the surface 
$c_{\varepsilon}(S_u)$, the process encounters the environment that is going to take it  through the polytope  $Q_u$ in time less than $4L/{\delta_0}$. 
Such an environment can be constructed in the following way. We require that the vector field $b$ satisfies $b\cdot l>\frac{\delta_0}2$ in the central cells of $Q_u$. The outer cells serve as traps that bring the process towards the central cells. More precisely if $C$ is one of the outer cells, we choose a vector $v_C$ that points from the center of $C$ towards the simplex that is shared with one of the interior cubes and that satisfies $v_C\cdot l>0$. 
We require $b$ to satisfy $$b\cdot v_C\geq \frac{\sqrt{d-1}}{\sqrt d}\cdot |b|\cdot |v_C|
\mbox{  
and }b\cdot v_C>\frac{\delta_0}2.$$

The environment now acts as a tunnel through which the process must go. This reduces the portion of the environment that gets exposed to the process, and plays the role of the sequence of steps of size $w$ that we used in proving (\ref{sub_add}).
In the end of the tunnel $Q_u$ we require that the environment has a set of teleports in the direction $u_1$, while before the end there are no teleports overlapping with the tube. 
More precisely, every point of the boundary of $Q_u$ that belongs to $Z_y$ is within $r/2$ of a $u_1$ teleport and at a distance of at least $r$ from all other teleports. 
Let $\hat T'_u$ be the time by which the process has spent time $t_0$ in the union of $u_1$ teleports. We assume that the value of the corresponding term $Y_{\lfloor \hat T'_u/t_0\rfloor}$ is not $\phi$. This means that $\hat T'_u$ is the time of jump. The conditional probabilities of such environments on $Q_u$ are positive, and we are sure that after the passage through the tunnel the process will appear at a uniform location within a cube that is sufficiently far away from the old environment. 

Let $B_2$ be the event that there exists $i\in\{1,2,\dots, m\}$ and $t_1\in\left(t_0\lfloor \frac{\hat T_u}{t_0}\rfloor,\hat T_u\right)$ such that $X_t$ spent all the time $(t_1-t_0, t_1)$ in a $u_i$ teleport. The definition of $B_2$ is the same as the definition of $
\hat F$ from the proof of Lemma \ref{first_lemma_cont}.

 We have the following sequence of inequalities
\begin{align*}
\mathbb P(T_{u+v}\leq &p+q, T_{u+v}\leq D_1)
\geq\mathbb
P\left(T_{u+v}\leq p+q, T_{u+v}\leq D_1, T_u\leq p\right)\\
&\geq\mathbb E\left[1(\hat T_u\leq p)\cdot 1(B_1)\cdot 1(T_{u+v}\leq p+q)\cdot 1(T_{u+v}\leq D_1)\right]
\\&=\mathbb E\Big[\mathbb E\Big[1(\hat T_u\leq p) \cdot 1(B_1) \cdot 1(T_{u+v}\leq p+q)\cdot 1(T_{u+v}\leq D_1)\big|\mathcal G_{\left[\hat T_u,\hat T'_u\right]}\Big]\Big]
,\end{align*}
where $\mathcal G_{[\hat T_u,\hat T'_u]}$ is the $\sigma$-algebra determined by the environment in the $M$-neighborhood of the path of the process from $0$ to $T_{u+v}$ except the the portion between times $\hat T_u$ and $\hat T'_u$. Let us introduce the following notation $\tilde X_t=X_{\hat T'_{u}+t}$. Let $\vec o$ be the vector that determines the translation which maps $[0,1]^d$ to the cube that contains $\tilde X_{0}$. Let $\tilde O$ be the translation of the origin for the vector $\vec o$. This way we can understand $\tilde X_t$ as a new process, and $\tilde O$ plays the role of the origin for that process. Let $\tilde Z_0$ be the hyperplane orthogonal to $l$ that passes through $\tilde O$ and define $\tilde D_1$ as the first time $\tilde X_t$ backtracks over $\tilde Z_0$. Define $\tilde T_v=\inf\{t: \tilde X_t\cdot l\geq v\}$.
Let $\mu=\frac{2L}{\delta_0}$. We know that $\hat T'_u-\hat T_u\leq \mu$ hence 
\begin{align*}
\mathbb P(T_{u+v}\leq p+q, &T_{u+v}\leq D_1)
\\
&\geq\mathbb E
\Big[\mathbb E\Big[1(\hat T_{u}\leq p) \cdot 1(B_1)
\cdot 1(\tilde T_{v}\leq q-\mu)\cdot 1(\tilde T_{v}\leq \tilde D_1)\big| 
\mathcal G_{\left[\hat T_u,\hat T'_u\right]}\Big]\Big].
\end{align*}
Since $\{\hat T_u\leq p\}$, $\{\tilde T_v\leq q-\mu\}$, $\{\tilde T_v\leq \tilde D_1\}$, $B_2$ $\in\mathcal G_{\left[\hat T_u,\hat T'_u\right]}$ we have
\begin{align*}
\mathbb E
\Big[\mathbb E\Big[&1(\hat T_{u}\leq p) \cdot 1(B_1)
 \cdot 1(\tilde T_{v}\leq q-\mu)\cdot 1(\tilde T_{v}\leq \tilde D_1)\big|\mathcal G_{\left[\hat T_u,\hat T'_u\right]}\Big]\Big]
\\
 &=
 \mathbb E
\Big[1(\hat T_u\leq p) 
 \cdot 1(\tilde T_{v}\leq q-\mu)\cdot 1(\hat T_u\leq D_1)
 \cdot 1(\tilde T_{v}\leq \tilde D_1)\cdot 1(B_2)\cdot
 \mathbb E\Big[1(B_1)
 \big|\mathcal G_{\left[\hat T_u,\hat T'_u\right]}\Big]\Big]\\
 &\;\;\; +
 \mathbb E
\Big[1(\hat T_u\leq p) 
 \cdot 1(\tilde T_{v}\leq q-\mu)\cdot 1(\hat T_u\leq D_1)\cdot 1(\tilde T_{v}\leq \tilde D_1)\cdot 1(B_2^C)\cdot
 \mathbb E\Big[1(B_1)
 \big|\mathcal G_{\left[\hat T_u,\hat T'_u\right]}\Big]\Big].
 \end{align*}
Now we bound each of the two terms on the right hand side in the same way as we did in the proof of Lemma \ref{first_lemma_cont} to get 
$$\mathbb P(T_{u+v}\leq p+q, T_{u+v}\leq D_1)\geq 
c\cdot \mathbb E(1(\hat T_u\leq p)\cdot 1(\hat T_u\leq D_1)\cdot 1(\tilde T_v\leq q)
\cdot 1(\tilde T_v\leq \tilde D_1)). $$

We can now replace $\tilde X$ with a process on the independent environment because it is sufficiently far away from $X$ and use the independence to obtain 
\begin{eqnarray*}
\mathbb P\left(T_{u+v}\leq p+q, T_{u+v}\leq D_1\right)&\geq&c\cdot
\mathbb P(\hat T_{u}\leq p, \hat T_u\leq D_1)\cdot \mathbb P(T_v\leq q-\mu, T_v\leq D_1)\\
&\geq&c'\cdot\mathbb P(T_{ u}\leq p, T_u\leq D_1)\cdot \mathbb P(T_v\leq q-\mu, T_v\leq D_1).
\end{eqnarray*}
We can now proceed in the same way as in the case of the deterministic walk and establish the existence of the limit in  (\ref{cont_probabilities}). The proof of the concavity of $\phi^l(k)$ is the same as in the discrete case. 

We will now prove that 
\begin{eqnarray}\label{moment_generating_lower}
\liminf \frac1{t}\log\mathbb E\left[ e^{\lambda \cdot X_t}\right]&\geq&
\max\left\{0,\sup_{k>0}\left\{k+\Phi\left(\frac{\lambda}{k}\right)\right\}\right\}.
\end{eqnarray}
Denote by $\Lambda(\lambda)$ the right-hand side of (\ref{moment_generating_lower}).
For each $k>0$ we have $\mathbb E[\exp\{\lambda\cdot X_t\}]\geq 
e^{tk}\mathbb P(X_t\cdot \lambda\geq kt)$ hence we only need to show that
$$\liminf \frac1t\log\mathbb E\left[e^{\lambda\cdot X_t}\right]\geq 0.$$
For this it suffices to show that $\frac1t\log\mathbb P\left(\lambda \cdot X_t\geq 0\right)\geq 0$ and for this we need to find a uniform lower bound on $\mathbb P(\lambda \cdot X_t\geq 0)$. 
In the beginning of the proof of Lemma \ref{first_lemma_cont} we showed that there is a positive probability that a given cube of size $3$ is a trap. Once the process enters there, it can exit only by teleportation. Hence there is a positive probability that the cube $[0,3]^d$ is a trap. Let $D$ be that event, and   
$D'$ the event that there are no teleports in proximity of $[0,3]^d$. Then we have $\mathbb P(D\cap D')\geq c$ for some positive constant $c$. The proof of (\ref{moment_generating_lower}) is complete once we observe that $\{\lambda\cdot X_t\geq 0\}\supseteq D\cap D'$. 

Proving inequality $$\limsup \frac1t\log\mathbb E\left[e^{\lambda\cdot X_t}\right]
\leq\max\left\{0,\sup_{k>0}\left\{k+\Phi\left(\frac{\lambda}k\right)\right\}\right\}$$ is exactly the same as in the discrete case. The proof only uses the concavity of $\Phi$ in the variable $k$. The convexity of $\Lambda$ is established also in the same way as in the discrete case, and this completes the proof of Theorem \ref{det_process}.\hfill $\Box$

\vspace{0.5cm}

We can use similar  arguments as in the discrete case to prove the following theorem: 

\begin{theorem}\label{mainmain_theorem_cont} Let $X_t$ be the previously defined 
process. 
Let $\Lambda$ be the function defined after (\ref{moment_generating_lower}) and let $\Lambda^*$ be 
its convex conjugate. Let $\mathcal F$ be the set of exposed points of $\Lambda^*$.
For any closed set $F\subseteq \mathbb R^d$, 
$$\limsup \frac1n \log\mathbb P\left(\frac1nX_n\in F\right)\leq -\inf_{x\in F} 
\Lambda^*(x),$$ 
and for any open set $G\subseteq \mathbb R^d$
$$\liminf \frac1n\log\mathbb P\left(\frac1nX_n\in G\right)\geq 
-\inf_{x\in G\cap \mathcal F} \Lambda^*(x).$$
Moreover, there exists $\delta>0$ such that $\Lambda^*(x)<+\infty$ for $|x|<\delta$. 
\end{theorem}
\noindent{\bf Proof.} 
The only difference from the proof in the discrete setting is proving that 
$\Lambda^*(x)<+\infty$ in some neighborhood of $0$. We follow the same idea, 
though. Denote by $e_i$ the base vectors of $\mathbb R^d$. Define $$\rho(l)=
\max\{\max_i\{e_i\cdot l\}, \max_i\{-e_i\cdot l\}\}.$$ Then $\inf_{l:|l|=1}\rho(l)>\frac1{\sqrt d}$. Hence, given a vector $\lambda$, we may assume that  $e_1\cdot \lambda> \frac1{\sqrt d}|\lambda|$. 
Consider the environment such that:
\begin{enumerate}
\item[(i)] The vector fields $b_Q$ satisfy $$b_Q\cdot e_1>\delta_0\;\mbox{ and }\;b_Q\cdot e_1>
\frac{\sqrt{d-1}}{\sqrt d}|b_Q|\;\;\mbox{ in the set }\;c_{\varepsilon}(Q),$$ for all unit cubes $Q$ whose lower-left corner has the form $(n,0,\dots, 0)$ where $n$ $\in$$\{0$, $1$, $\dots$, $\lfloor t \delta_0\rfloor\}$. 
\item[(ii)] In each unit cube $Q$ that shares a boundary to any of the cubes from (i) 
we consider the vector $v\in\{-1,0,1\}^d\setminus\{0\}$ such that 
$v\cdot \lambda$ is maximal. Then we require $b_Q$ to satisfy $$b_Q\cdot v\geq \delta_0|v|\;\mbox{ and }\;
b_Q\cdot v\geq \frac{\sqrt{d-1}}{\sqrt d}|b_Q||v|.$$
\item[(iii)]  Denote by  $Q_t$ the cube whose lower-left corner is $(\lfloor t\delta_0\rfloor +1,0,0,\dots,0)$. We require that cells adjacent to $Q_t$ serve as a trap towards $Q_t$. This environment is precisely described in the beginning of the proof of Lemma \ref{first_lemma_cont}.
\item[(iv)] There are no teleports in the $5\sqrt d +r$ neighborhood of the segment 
from $O$ to the lower-left corner of $Q_t$.
\end{enumerate}

The conditional probabilities of the previously described environment are bounded 
below by $c^{t}$ for some $c>0$. On the realization of such environment we can 
guarantee that $$X_t\cdot \lambda>
\frac{t\delta_0|\lambda|}{2\sqrt d}.$$ Hence 
\begin{eqnarray*}\Phi\left(\frac{\lambda}{\delta_0|\lambda|\cdot\frac1{2\sqrt d}}\right)&=&\lim_{t\to+\infty}\frac1t\log\mathbb P\left(X_t\cdot \lambda> \frac{t\delta_0|\lambda|}{2\sqrt d}\right)\geq c.\end{eqnarray*}
Therefore $\Lambda(\lambda)\geq \frac{\delta_0|\lambda|}{2\sqrt d}+c$ and 
$$\Lambda^*(x)\leq \sup_x\left\{\lambda\cdot x-\frac{\delta_0|\lambda|}{2\sqrt d}-
c\right\}\leq
\sup_x\left\{|\lambda|\cdot |x|-\frac{\delta_0|\lambda|}{2\sqrt d}\right\}-c.
$$
If $|x|<\frac{\delta_0}{2\sqrt d}$ we see that $\Lambda^*(x)\leq-c$.
\hfill$\Box$

\vspace{0.5cm}

\subsection{Open Problems} One of the main unsolved questions is the continuity of $\Phi$. We are only able to prove that $k\mapsto \Phi\left(\frac lk\right)$ is concave for fixed $l\in\mathbb R^d$, which gives the continuity in $k$. However, it would be very interesting to see the continuity of the $\Phi$ as a function defined on a domain from $\mathbb R^d$. Ideas developed in this paper can't be directly applied because  the hitting times of hyperplanes can't be easily replaced by the hitting times of dihedral angles (or some other interfaces) in the proofs of the previous theorems.  The said difficulty appears even in the discrete case. 

\vspace{0.5cm}
\noindent{\em Acknowledgements.} I would like to thank my advisor, professor Fraydoun Rezakhanlou, for formulating the probability model that was studied here and for  excellent ideas and advices on how to approach the problem. I would also like to thank  Firas Rassoul-Agha and Atilla Yilmaz for numerous discussions.
I am also very grateful to the referee for pointing out mistakes and for giving nice comments and ideas for simplifying certain proofs.

\end{document}